\documentstyle{amsppt}
\topmatter
\title Companion forms and the structure of $p$-adic Hecke algebras
\endtitle
\author Masami Ohta \endauthor
\endtopmatter
\document
\vskip5mm
{\bf Abstract.} We study the structure of the Eiesenstein
component of Hida's ordinary $p$-adic Hecke algebra
attached to {\it modular forms}, in connection with the 
companion forms in the space of modular forms (mod $p$).
We show that such an algebra is a Gorenstein ring
if certain space of modular forms (mod $p$) having
companions is one-dimensional; and also give a numerical
criterion for this one-dimensionality. This in part overlaps
with a work of Skinner and Wiles; but our method,
based on a work of Ulmer, is totally different.
We then consider consequences of the above mentioned Gorenstein property.
We especially discuss the connection with the Iwasawa theory.

\vskip3mm
{\bf Mathematics Subject Classification 2000.} 11F33, 11F80.

\vskip5mm
{\bf Introduction}

\vskip3mm
The connection of the modular Galois representations with the
theory of cyclotomic fields was first studied by Ribet [R].
Then the Iwasawa main conjecture for ${\bold Q}$ was proved by Mazur and 
Wiles in their fundamental paper [MW]. After these works, Harder and Pink
[HP] and Kurihara [Ku] have independently shown that, if the Eisenstein
component of certain $p$-adic Hecke algebra attached to {\it cusp forms}
is a Gorenstein ring, then one can further determine the
structure of certain Iwasawa module associated with ${\bold Q}
(\mu_{p^{\infty}})/{\bold Q}$ completely.

Although the Gorenstein property 
is established for most of non-Eisenstein components 
of such $p$-adic Hecke algebras by the works of
many mathematicians (cf. Wiles [Wi], Chapter 2, Section 1), 
little is known about the Eisenstein components.

In contrast to this, the same property of an Eisenstein component of
the $p$-adic Hecke algebra attached to {\it modular forms} seems
relatively easier to achieve. Indeed, in their course of the
proof of 
$`` R \cong {\bold T}"$-type theorem for residually reducible
two-dimensional Galois representations, Skinner and Wiles [SW]
have shown that such an algebra is even a complete intersection,
under some hypotheses. 

The aim of the present article is to investigate this latter problem by
a method totally different from [SW], and then to consider some
consequences. 

\vskip3mm
Though we treat the
Hecke algebras of level $N$ with $p \nmid \varphi (N)$ in the text,
we would like to explain the contents of this paper under the 
assumption $N=1$ for simplicity, adding comments on the general case.

We throughout fix a prime number $p \geq 5$. Let $e\, M$
(resp. $e\, S 
)$ be the space of ordinary $\Lambda$-adic modular forms
(resp. $\Lambda$-adic cusp forms), introduced by Hida and Wiles,
of level 1 over ${\bold Z}_p$. (This is the space denoted by
$e\, M(1; \Lambda_{{\bold Z}_p})$
(resp. $e\, S(1;\Lambda_{{\bold Z}_p})$) in the text.)
Let $e\, {\Cal H}=e\, {\Cal H}(1;{\bold Z}_p)$ 
(resp. $e\, h =e\, h(1; {\bold Z}_p)$) be
Hida's universal ordinary $p$-adic Hecke algebra acting on
this space. It is a fundamental fact, due to Hida, that these spaces 
and algebras are finite and flat over the Iwasawa algebra
$\Lambda ={\bold Z}_p[[T]]$.
They split into a direct sum of
$\omega^i$-eigenspaces with respect to the action of 
$({\bold Z}/p{\bold Z})^{\times}$, where
$\omega$ denotes the Teichm\"uller character. We indicate 
these eigenspaces by the superscript $``{}^{(i)}"$.

Fix an even integer $k$ such that $2 \leq k 
\leq p-1$. The space $e\, M^{(k-2)}$ contains an essentially unique
$\Lambda$-adic Eisenstein series 
which interpolates the classical Eisenstein series.
We define the {\it Eisenstein ideal} 
${\Cal I}$ of $e\, {\Cal H}^{(k-2)}$ as 
the annihilator of this $\Lambda$-adic Eisenstein
series.
Then ${\frak M}:=
({\Cal I}, p ,T)$ is a maximal ideal of $e\, {\Cal H}^{(k-2)}$,
and our main concern is the localization $e\, {\Cal H}_{\frak M}^{(k-2)}$,
the {\it Eisenstein component} of $e\, {\Cal H}^{(k-2)}$. 
The Eisenstein component of $e\, h^{(k-2)}$ is then defined
as $e\, h_{\frak M}^{(k-2)}:= e\, h^{(k-2)}
\otimes_{e\, {\Cal H}^{(k-2)}}e\, {\Cal H}_{\frak M}^{(k-2)}$.
(These are the algebras denoted by 
$e\, {\Cal H}(1; {\bold Z}_p)_{{\frak M}(\omega^{k-2}, {\pmb 1})}$ 
and $e\,h(1;{\bold Z}_p)_{{\frak M}(\omega^{k-2},{\pmb 1})}$
in the text, respectively.)
This latter ring is not a zero-ring if and only if the Bernoulli
number $B_k$ is not a $p$-adic unit.
We have $e\, {\Cal H}_{\frak M}^{(k-2)}=\Lambda$ 
if $k=2$ or $p-1$,
and hence we assume that $4 \leq k\leq p-3$, in the rest of this introduction.

\vskip1mm
Let $\widetilde{M}_k$ (resp. $\widetilde{S}_k$) $\subset {\bold F}_p
[[q]]$ be the space of modular forms (resp. cusp forms) (mod $p$) 
of weight $k$ and level 1 in the sense of Serre and 
Swinnerton-Dyer. We can consider the Hecke algebra
acting on $\widetilde{M}_k$, i.e. the subalgebra of
${\roman{End}}_{{\bold F}_p}(\widetilde{M}_k)$ generated 
by all Hecke operators
$T(n)$. Let $\overline{\frak m}$ be the annihilator of the unique Eisenstein
series (mod $p$) in $\widetilde{M}_k$, 
in this algebra.
The Eisenstein component 
of $\widetilde{M}_k$ is then
the localization
$\widetilde{M}_{k,\overline{\frak m}}$. 
Set $k':=p+1-k$, and define $\overline{\frak m}'$ 
and $\widetilde{M}_{k', \overline{\frak m}'}$
from $\widetilde{M}_{k'}$ in the same manner as above.
(The spaces $\widetilde{M}_{k, \overline{\frak m}}$ and
$\widetilde{M}_{k',\overline{\frak m}'}$
are the ones denoted by $M_k(\varGamma_1(1);{\bold F}_p)_{\frak m}$
and $M_{k'}(\varGamma_1(1);{\bold F}_p)_{{\frak m}'}$ in the text,
respectively. The ideals $\overline{\frak m}$ and $\overline{\frak m}'$
are the images of ${\frak m}$ and ${\frak m}'$ in the text in the
characteristic $p$ Hecke algebra.)

Now recall that $f\in \widetilde{M}_k$ and
$g \in \widetilde{M}_{k'}$ are said to be companions to each other if
$$\theta^{k'}f=\theta g \quad {\roman{or}}\enskip{\roman{equivalently}}
\quad \theta^kg=\theta f$$
where $\theta :=q(d/dq)$ is the usual derivation.
Let $c(\overline{\frak m}')$ (which will be denoted by $c({\frak m}')$
in the text)
be the dimension of 
the space $\{ g \in \widetilde{M}_{k',\overline{\frak m}'}
\mid {\text{$g$ has a companion}}\}$ over ${\bold F}_p$.
It is easy to see that $c(\overline{\frak m}')$ is always positive.

\proclaim{Theorem 1} If $c(\overline{\frak m}')=1$, then
$e\, {\Cal H}_{\frak M}^{(k-2)}$ is a Gorenstein ring.
\endproclaim

The nature of the space above is not clear in general
at present, but here is a simple criterion:
$$B_{k'}\in {\bold Z}_p^{\times}
\quad \Rightarrow \quad c(\overline{\frak m}')=1. \tag $*$ $$
Therefore the left numerical condition implies that 
$e\, {\Cal H}_{\frak M}^{(k-2)}$ is Gorenstein.
We note that the above quoted work of Skinner and Wiles 
shows that the
same condition implies that $e\,{\Cal H}_{\frak M}^{(k-2)}$ 
is moreover a complete
intersection. In the text, however, we give a result as in Theorem 1
for a wider class of Eisenstein components than treated in [SW]. The 
numerical criterion $(*)$ is also generalized,
using our previous work [O4].
 
Our proof of Theorem 1 is based on the following result due to
Ulmer [U]:
\proclaim{Theorem} There is a bilinear pairing
$$\widetilde{M}_k \times \widetilde{S}_k\to {\bold F}_p$$
such that:

{\rm i)} its left kernel consists of the forms with companions;

{\rm ii)} all Hecke operators are self-adjoint with respect
to it.
\endproclaim

Ulmer constructed his pairing for general level $N$,
in the course of detailed study of the Selmer group
of the universal elliptic curve over the Igusa curve. However, we
remark that the original pairing, as it stands, does not satisfy
ii) above when $N>1$; but it is not difficult to get a Hecke equivariant
pairing by twisting it.

One of the motivation of this work was to study the connection of the objects
above with the theory of cyclotomic fields. 
Let $L_{\infty}$ be the maximal abelian unramified pro-$p$ extension
of ${\bold Q}(\mu_{p^{\infty}})$. As before, we denote by
${\roman{Gal}}(L_{\infty}/{\bold Q}(\mu_{p^{\infty}}))^{(i)}$
the $\omega^i$-eigenspace of the associated 
Galois group with respect to the action of
${\roman{Gal}}({\bold Q}(\mu_p)/{\bold Q})$.
If $e\, h_{\frak M}^{(k-2)}$ is a zero-ring, 
$e\, {\Cal H}_{\frak M}^{(k-2)}=\Lambda$, ${\Cal I}$ is its
principal ideal, 
and ${\roman{Gal}}(L_{\infty}/{\bold Q}(\mu_{p^{\infty}}))^{(1-k)}
=\{ 0 \}$. We henceforth assume otherwise.

\proclaim{Theorem 2} If $e\, {\Cal H}_{\frak M}^{(k-2)}$ 
is a Gorenstein ring, then the following conditions are equivalent:

{\rm i)} $e\, h_{\frak M}^{(k-2)}$ is Gorenstein;

{\rm ii)} $e\, h_{\frak M}^{(k-2)}$ is a complete intersection;

{\rm iii)} the Eisenstein ideal ${\Cal I}$ of $e\, {\Cal H}^{(k-2)}$ is
principal;

{\rm iv)} the image of ${\Cal I}$ in $e\, h^{(k-2)}$ is principal;

{\rm v)} the Iwasawa module ${\roman{Gal}}
(L_{\infty}/{\bold Q}(\mu_{p^{\infty}}))^{(1-k)}$ is cyclic over $\Lambda$.
\endproclaim
The implications: iii) $\Leftrightarrow$ iv) $\Rightarrow$ ii)
$\Rightarrow$ i) are easy or well-known, and Harder, Pink and Kurihara proved
that i) $\Rightarrow$ v). These hold unconditionally.
The new point here is
that, under the hypothesis stated in the theorem, v) in turn implies
iv). More precisely, under the same hypothesis, we 
{\it explicitly describe} the
$\Lambda$-module ${\roman{Gal}}(L_{\infty}/{\bold Q}
(\mu_{p^{\infty}}))^{(1-k)}$ in terms of the Eisenstein ideal,
by applying the method of [HP] and [Ku]
to the modular Galois representation 
into $GL_2(e\, h_{\frak M}^{(k-2)})$ attached to the Eisenstein
component of $e\, S^{(k-2)}$.

\vskip3mm
The outline of the text is as follows:

In Section 1, 
after notational preliminaries, we consider the comparison between
the space of ordinary modular forms of level $N$ and that of level $Np$;
and the duality between the space of modular forms and its 
Hecke algebra, for later use.
As for the latter, our main result is Proposition (1.5.3). 
It especially generalizes the known duality 
for the Eisenstein
components of the space of ordinary $\Lambda$-adic modular forms
([O4]) by a much simpler argument.

In Section 2, we study Ulmer's pairing and its twisted version. We
first recall the definition of the pairing given in [U], and then look at the
behavior of the Hecke operators and the diamond operators 
with respect to it.  
For general level $N$, we show that the pairing obtained 
by twisting Ulmer's original one by an operator $w_N$ is
Hecke equivariant. It is stated as
Theorem (2.5.1) together with its complement (2.5.5).
This procedure is quite analogous to the 
construction of the Hecke equivariant
$``$twisted Weil pairing" out of the Weil pairing.

With these two preliminary sections, we are ready to prove our
main results in Section 3. The general forms of Theorem 1 and
the numerical criterion
$(*)$ are given by Theorems (3.3.2) and (3.2.3), respectively.
We also show, in Theorem (3.3.8),
that the implication: i) $\Rightarrow$ iii) as in
Theorem 2 follows from the hypothesis as in Theorem 1.
We finally consider the connection with the theory of cyclotomic
fields, the general form of Theorem 2 being presented as
Theorem (3.4.12) and Corollary (3.4.13).

\vskip5mm
\flushpar
{\bf \S1. Preliminaries on modular forms}

\vskip3mm \flushpar
{\bf 1.1. Notational preliminaries.}
Let $k$ be an integer such that $k \geq 2$.
We denote by $M_k(\varGamma_1(M))$
(resp. $S_k(\varGamma_1(M))$) the 
complex vector space of modular forms
(resp. cusp forms) of weight $k$
with respect to $\varGamma_1(M)$. Via the $q$-expansion
at infinity, we consider it as a subspace of ${\bold C}[[q]]$, and
set

$$ \cases
M_k(\varGamma_1(M); {\bold Z}):= M_k(\varGamma_1(M))\cap {\bold Z}[[q]] \\
S_k(\varGamma_1(M);{\bold Z}):=S_k(\varGamma_1(M))\cap {\bold Z}[[q]].
\endcases \tag 1.1.1$$
Also, for any ring $R$, we set
$$\cases
M_k(\varGamma_1(M); R):= M_k(\varGamma_1(M);{\bold Z})\otimes_{\bold Z} R 
\hookrightarrow R[[q]]\\
S_k(\varGamma_1(M);R):= S_k(\varGamma_1(M);{\bold Z})\otimes_{\bold Z} R
\hookrightarrow R[[q]]
\endcases \tag 1.1.2 $$
so that the base changing property

$$\cases
M_k(\varGamma_1(M); R')=M_k(\varGamma_1(M); R)\otimes_R R' \\
S_k(\varGamma_1(M);R')=S_k(\varGamma_1(M); R)\otimes_R R'
\endcases \tag 1.1.3$$
trivially holds for any ring homomorphism $R\to R'$.

When $R$ is a ${\bold Z}[1/6M]$-algebra,
the spaces in (1.1.2) are canonically isomorphic to the spaces of
modular forms or cusp forms
defined geometrically as in Katz [Ka1] (cf. also Gross [Gr]), via the
$q$-expansion mapping.
This is a consequence of
the $q$-expansion principle and the base changing property
as (1.1.3) for such forms
(cf. [Gr], Proposition 2.5 and remarks in Section 10).
Later, when it is convenient to
do so, we will take this geometric point of view. For the
argument of this section, however,
the above formal definition is rather convenient.
If $R$ is a finite field of
characteristic $p$, these are the spaces of modular forms
or cusp forms (mod $p$) over $R$ in the sense of Serre
and Swinnerton-Dyer.

The usual Hecke operators $T(n)$ (for positive integers $n$)
and $T(m,m)$ (for positive integers $m$ prime to $M$), as well as the
diamond operators $\langle m \rangle$ (for integers $m$ 
prime to $M$) preserve $M_k(\varGamma_1(M);{\bold Z})$
and $S_k(\varGamma_1(M); {\bold Z})$. We then define the Hecke algebras

$$\cases H_k(\varGamma_1(M); R) \\
h_k(\varGamma_1(M); R) \endcases \tag 1.1.4 $$
as the $R$-subalgebra of ${\text {End}}_R(M_k(\varGamma_1(M); R))$ or
${\text {End}}_R(S_k(\varGamma_1(M); R))$ generated by all 
such operators. Here we recall that $T(m,m)=m^{k-2} \langle m \rangle$
for any positive integer $m$ prime to $M$.

Similar terminologies as above will be used for other congruence
subgroups of $SL_2({\bold Z})$.

\vskip3mm
In this paper, {\it we throughout fix a prime number $p \geq 5$
and a positive integer $N$ prime to $p$}.
Let $F$ be a finite extension of ${\bold Q}_p$ with its ring of integers
${\frak r}$. We denote by $\varpi$ a prime element of $F$, and by 
${\frak k}:={\frak r}/(\varpi )$ the 
residue field. When $R=F, {\frak r}$ or ${\frak k}$,
the Hecke algebras above are generated over $R$ by all $T(n)$; or by all
$T(l)$ with prime numbers $l$ and all $T(m,m)$.

Let $\Lambda_{\frak r}$ be the Iwasawa algebra
over ${\frak r}$, i.e. the completed group algebra over ${\frak r}$
of the multiplicative group $1+p{\bold Z}_p$. Fixing a topological
generator $\gamma$ of $1+p{\bold Z}_p$, we identify it with
the formal power series ring ${\frak r}[[T]]$ via $\gamma \leftrightarrow
1+T$.
Now we denote by 
$$\cases e\, M(N; \Lambda_{\frak r}) \\
e\, S(N; \Lambda_{\frak r}) \endcases \tag 1.1.5$$
the spaces of ordinary $\Lambda_{\frak r}$-adic modular forms and cusp forms 
introduced by Hida and Wiles. We use the same convention for such
spaces as in our previous works, and refer the reader to [O1] and [O2]
for details. 
We also refer the reader to
these papers for the definitions and basic properties
of Hida's universal ordinary $p$-adic Hecke algebras:
$$\cases e\, {\Cal H}(N; {\frak r}) \\
e\, h(N;{\frak r}) \endcases \tag 1.1.6$$
acting on the spaces in (1.1.5).
We only remind of us that they are algebras over
${\frak r}[[({\bold Z}/Np{\bold Z})^{\times}\times (1+p{\bold Z}_p)]]$
in such a way that a positive integer $m$ prime to $Np$, considered
as an element of $({\bold Z}/Np{\bold Z})^{\times}\times (1+p{\bold Z}_p)$,
acts as multiplication by $T(m,m)$.

\vskip3mm \flushpar
{\bf 1.2. Eisenstein series.} For Dirichlet characters $\chi$,
$\psi$ and a positive integer $c$, we formally set
$$E_k(\chi , \psi; c):= \delta(\psi )L(1-k, \chi)+
\sum_{n=1}^{\infty}\left( \sum_{0<t \mid n}\chi (t)\psi (\frac{n}{t})
t^{k-1}\right) q^{cn}. \tag 1.2.1$$
Here, $\delta (\psi )= 1/2$ or 0 according as $\psi$ is the trivial
character $\pmb{1}$ or not, and $L(s, \chi)$ is the Dirichlet $L$-function.
Let $\chi$ and $\psi$ be primitive 
of conductors $u_0$ and $v$, respectively.
Then it is well-known that $E_k(\chi ,\psi ;c)$ belongs 
to $M_k(\varGamma_1(N))$
whenever the following conditions:
$$\cases
(\chi \psi )(-1)= (-1)^k; \\
{\text{$cu_0v$ divides $N$}}
\endcases \tag 1.2.2 $$
are satisfied, unless $k=2$ and $\chi = \psi ={\pmb{1}}$.

Next we recall $\Lambda$-adic Eisenstein series (cf. [O2], 2.3).
Let $\psi$,
$v$ and $c$ be
as above and let $\vartheta$ be a primitive Dirichlet character of
conductor $u$. We assume that $F$ contains the values of $\vartheta$
and $\psi$. Then the series
$${\Cal E}(\vartheta , \psi ;c) :=
\delta (\psi )G(T,\vartheta \omega^2)+
\sum_{n=1}^{\infty}\left( \underset{p\nmid t}\to{\sum_{0 <t\mid n}}
\vartheta (t)\psi (\frac{n}{t})A_t(T)\right )q^{cn} \in \Lambda_{\frak r}[[q]]
\tag 1.2.3 $$
belongs to $e\, M(N;\Lambda_{\frak r})$ whenever

$$\cases
\vartheta \psi (-1)=1;\\
{\text{$v$ and $c$ are prime to $p$}};\\
{\text{$cuv$ divides $Np$;}} \\
(\vartheta , \psi )\not= (\omega^{-2}, {\pmb 1}).
\endcases \tag 1.2.4$$
Here and henceforth, 
$\omega$ denotes the Teichm\"uller character,
$$A_t(T):=t(1+T)^{s(t)}~~\text{ if }~~t\omega^{-1}(t)=\gamma^{s(t)}
\tag 1.2.5$$
and $G(T;\vartheta \omega^{2})\in \Lambda_{\frak r}$ 
is a twist of Iwasawa's power series
giving the Kubota-Leopoldt $p$-adic $L$-function:
$$G(\gamma^s-1 ,\vartheta \omega^2)=L_p(-1-s, \vartheta \omega^2).
\tag 1.2.6$$

For a character
$\varepsilon$ of $1+p{\bold Z}_p$ of finite
order and an integer $d \geq 0$, the $(\varepsilon ,d)$-specialization
of ${\Cal E}(\vartheta ,\psi ;c)$,
i.e. the classical form obtained by setting $T=\varepsilon (\gamma )
\gamma^d-1$, is given by:
$${\Cal E}(\vartheta ,\psi ;c)_{\varepsilon , d}=
E_{d+2}((\vartheta \varepsilon \omega^{-d})_1, \psi ;c) 
\tag 1.2.7 $$
where the subscript $``{}_1"$ indicates the (possibly imprimitive) 
induced character
defined modulo the least common multiple of its conductor and $p$.
If $\varepsilon$ takes values in $F$ and
Ker$(\varepsilon )=1+p^r{\bold Z}_p$, this series belongs to
$M_{d+2}(\varGamma_1(Np^r); {\frak r})$.
When $(\vartheta ,\psi )=(\omega^{-2}, {\pmb 1})$, the form
$((1+T)-\gamma^{-2}){\Cal E}(\omega^{-2}, {\pmb 1};c)\in
e\, M(N;\Lambda_{\frak r})$ has a similar property.

Now assume $uv=N$ or $Np$, and set
$${\Cal E}(\vartheta ,\psi ):={\Cal E}(\vartheta ,\psi ; 1).
\tag 1.2.8$$
It is a common eigen form of all $T(n)$ and $T(m,m)$ in 
$e\,{\Cal H}(N;{\frak r})$, and we let
$$\cases
{\Cal I}(\vartheta ,\psi )= {\Cal I}:=
{\roman{Ann}}_{e\, {\Cal H}(N;{\frak r})}({\Cal E}(\vartheta ,\psi ))
,{\text{ the annihilator of ${\Cal E}(\vartheta ,\psi )$ in
$e\,{\Cal H}(N;{\frak r})$}};\\
{\frak M}(\vartheta ,\psi )={\frak M}:= ({\Cal I}(\vartheta ,\psi ),\varpi ,T)
\endcases \tag 1.2.9$$
be the {\it Eisenstein ideal} and the {\it Eisenstein maximal ideal}
of $e\, {\Cal H}(N;{\frak r})$
attached to ${\Cal E}(\vartheta ,\psi )$, respectively. 
The ideal ${\Cal I}(\vartheta ,\psi )$
is generated by all 
$$T(n) - 
\underset{p \nmid t}\to{\sum_{0<t\mid n}}\vartheta (t)\psi (n/t)
A_t(T).$$
For any $e\, {\Cal H}(N;{\frak r})$-module, we will indicate by 
the subscript $``{}_{\frak M}$" the localization at ${\frak M}$.
 
For a non-negative integer $d$, we set
$$\omega_d :=T-({\gamma}^d -1) \in \Lambda_{\frak r}.
\tag 1.2.10$$
Then there is a canonical isomorphism:
$$e\, {\Cal H}(N;{\frak r})/\omega_d \overset{\sim}\to\rightarrow
e\, H_{d+2}(\varGamma_1(Np);{\frak r})
\tag 1.2.11$$
$e$ in the right hand side being Hida's idempotent attached to $T(p)$.
The image of ${\frak M}(\vartheta ,\psi )$ in $e\, H_{d+2}(\varGamma_1(Np);
{\frak r})$ under this isomorphism
is the Eisenstein maximal ideal defined similarly as above with respect
to $E_{d+2}((\vartheta \omega^{-d})_1, \psi ;1)$.
Thus, if $X$ is an $e\, H_{d+2}(\varGamma_1(Np);{\frak r})$-module,
$X_{\frak M}$ may be identified with its localization at this Eisenstein
maximal ideal.

\vskip3mm \flushpar
{\bf 1.3. Ordinary modular forms of level $N$ and 
level $Np$}. Clearly,
$S_k(\varGamma_1(N);{\frak r})$ (resp. $M_k(\varGamma_1(N);{\frak r}))$
is an ${\frak r}$-submodule of $S_k(\varGamma_1(Np);{\frak r})$ 
(resp. $M_k(\varGamma_1(Np); {\frak r})$). 
In the following, we denote by $e$ Hida's idempotent
attached to $T(p)$ on the latter space.
The first statement of the following lemma is due to 
Gouv{$\hat{\roman{e}}$}a
[Go], Proposition 5:

\proclaim{Lemma (1.3.1)} When $k \geq 3$, we have
$$\cases
e\, S_k(\varGamma_1(N); {\frak r})=
e\, S_k(\varGamma_1(N) \cap \varGamma_0(p);{\frak r})\\
e\, M_k(\varGamma_1(N);{\frak r})=
e\, M_k(\varGamma_1(N)\cap \varGamma_0(p); {\frak r}).
\endcases $$
\endproclaim

\demo{Proof} We deduce the second statement from Gouv$\hat{\roman e}$a's
result. For this, it is enough to show that $e\, M_k(\varGamma_1(N);
F)\supseteq e\, M_k(\varGamma_1(N)\cap \varGamma_0(p);
F)$ when $F$ contains all $\varphi (N)$-th roots of unity, $\varphi$ being the
Euler function.

It is then easy to see that $E_k(\chi_1, \psi ;c)$ with primitive $\chi$,
$\psi$ and $c$ satisfying (1.2.2), together with ordinary cusp
forms span $e\, M_k(\varGamma_1(N)\cap \varGamma_0(p); F)$ (cf. [O2], Lemma
(2.3.2), where a result of Hida is recalled). But it follows from
Hida [H1], Lemma 3.3 that $E_k(\chi_1, \psi ;c)$ is a non-zero
constant multiple of $e\, E_k(\chi ,\psi ;c)$, and hence belongs
to $e\, M_k(\varGamma_1(N);F)$.
\qed \enddemo  

Let $e_0$ be the idempotent attached to $T(p)\in H_k(\varGamma_1(N);
{\frak r}).$

\proclaim{Proposition (1.3.2)} Assume $k \geq 3$. Then $e$ induces an
isomorphism
$$e: e_0\,M_k(\varGamma_1(N);{\frak r}) \overset{\sim}\to\rightarrow
e\, M_k(\varGamma_1(N)\cap \varGamma_0(p);{\frak r}).$$
This isomorphism commutes with $T(l)$ for primes $l \not= p$ and
$\langle m \rangle $ with $(m, Np )=1$. Moreover, $T(p)$ on the
left commutes with $T(p)+p^{k-1}\langle n_p \rangle T(p)^{-1}$ on
the right. Here, $n_p$ is an integer prime to $Np$ and congruent
to $p$ modulo $N$.
\endproclaim
\demo{Proof}
Gouv$\hat{\roman{e}}$a has shown that
$$e: e_0\,S_k(\varGamma_1(N); {\frak r}) \to
e\, S_k(\varGamma_1(N); {\frak r})$$
is an isomorphism, and also that the compatibilities for the Hecke 
operators as above hold ([Go], Lemma 3). The same argument works for
modular forms; and our claim follows from
the previous lemma.
\qed \enddemo

In the following corollaries, we also assume that $k \geq 3$.
\proclaim{Corollary (1.3.3)} There is an isomorphism of
${\frak r}$-algebras 
$$e_0\, H_k(\varGamma_1(N); {\frak r}) \overset{\sim}\to\rightarrow
e\, H_k(\varGamma_1(N)\cap \varGamma_0(p); {\frak r})$$
such that:
$$\cases
\text{$T(l)\mapsto T(l)$ for primes $l \not= p$;} \\
\text{$\langle m \rangle \mapsto \langle m \rangle$ for $m$ prime to
$Np$;}\\
T(p) \mapsto T(p)+p^{k-1}\langle n_p \rangle T(p)^{-1}.  \endcases  $$
\endproclaim

Since $T(p)$ of level $N$ and level $Np$ coincide in
characteristic $p$, we obtain the following:

\proclaim{Corollary (1.3.4)} We have the identity
in ${\frak k}[[q]]$:
$$e_0\, M_k(\varGamma_1(N); {\frak k})=
e\, M_k(\varGamma_1(N)\cap \varGamma_0(p); {\frak k}).$$
\endproclaim

\proclaim{Proposition (1.3.5)} For any integer $k \geq 3$, we have:
$$
e_0\, M_k(\varGamma_1(N); {\frak k})
= e\, M_2 (\varGamma_1(Np), \omega^{k-2};{\frak k}) $$
where the right hand side denotes the maximum direct summand
of 
$e\, M_2(\varGamma_1(Np);
{\frak k})$ on which the diamond operator
$\langle m \rangle$ acts as $\omega^{k-2} (m)$ for
$m \in \! ({\bold Z}/p{\bold Z})^{\times}\! 
\subseteq \! ({\bold Z}/Np{\bold Z})^{\times}$.
\endproclaim
\demo{Proof} Let $e\, M(N;\Lambda_{\frak r})^{(j)}$ be the 
$\omega^j$-eigenspace with respect to the action of 

\flushpar
$({\bold Z}/p{\bold Z})^{\times}$ (cf. 1.1).
It is then easy to see that
$({\pmb 1}, a)$-specialization gives an isomorphism:
$$e\, M(N;\Lambda_{\frak r})^{(k-2)}/\omega_a 
\overset{\sim}\to \rightarrow e\, M_{a+2}(\varGamma_1(Np),
\omega^{k-2-a};{\frak r})$$ 
for each integer $a \geq 0$
(cf. [O2], Proposition (2.5.4)), with the same meaning of the
right hand side as above.
Applying this to $a=0$ and $k-2$, and
further reducing modulo $\varpi$, we conclude that
$e\, M_k(\varGamma_1(N)\cap \varGamma_0(p); {\frak k})
=e\, M_2(\varGamma_1(Np), \omega^{k-2} ; {\frak k})$. Our result
follows from (1.3.4).
\qed \enddemo

For a related lifting property of not necessarily ordinary cusp
forms, cf. [Gr], Proposition 9.3. 
When $k=p+1$, the identity above reads as: 
$e_0\, M_{p+1}(\varGamma_1(N); {\frak k})=e\, 
M_2(\varGamma_1(N)\cap \varGamma_0(p) ;{\frak k}).$
This is a part of Serre [Se], Th\'eor\`eme 11, when $N=1$.

\vskip3mm \flushpar
{\bf 1.4. Duality between modular forms and Hecke algebras.} In the
following, we denote by $a(n; f)\in R$ the coefficient of
$q^n$ in $f \in M_k(\varGamma_1(N);R)\hookrightarrow R[[q]]$.
\proclaim{Proposition (1.4.1)} The pairing:
$$H_k(\varGamma_1(N);{\frak r}) \times
M_k(\varGamma_1(N);{\frak r})\to {\frak r}$$
defined by $(t , f):=a(1;f\mid t)$
gives a \, perfect pairing of free ${\frak r}$-modules, when
$2 \leq k \not\equiv 0$ $({\roman{mod}}$ $p-1)$.

Also the pairing:
$$e\,H_k(\varGamma_1(N)\cap \varGamma_0(p);{\frak r}) \times
e\, M_k(\varGamma_1(N)\cap \varGamma_0(p);{\frak r}) \to
{\frak r}$$
defined by the same formula
is perfect when $3 \leq k \not\equiv 0$ $(\roman{mod}$ $p-1)$.
\endproclaim
\demo{Proof} It is well-known that the formula above gives a 
perfect pairing
between $H_k(\varGamma_1(N);F)$ and
$M_k(\varGamma_1(N); F)$ (cf. Hida [H2], \S2). It therefore
induces an injection:
$$M_k(\varGamma_1(N);{\frak r}) \to {\roman{Hom}}_{\frak r}
(H_k(\varGamma_1(N);{\frak r}),{\frak r})$$
and what we need to show is its surjectivity. If $\phi$ is an
element of the right hand side, there is an $f\in M_k(\varGamma_1(N);
F)$ such that $\phi (t)= (t,f)$. We have $(T(n),f)=
a(n;f) \in {\frak r}$ for all $n \geq 1$. 
If $a(0;f)$ does not belong to ${\frak r}$, then
multiplying $f$ by a suitable power of $\varpi$
and reducing modulo $\varpi$, we see that a non-zero constant belongs
to $M_k(\varGamma_1(N);{\frak k})$, which can happen only when
$k$ is divisible by $p-1$.
This proves the first assertion.

The same proof works for the second assertion by virtue of (1.3.4).
\qed \enddemo

\proclaim{Corollary (1.4.2)} We have canonical ring isomorphisms:
$$\cases
H_k(\varGamma_1(N);{\frak r})/\varpi \overset{\sim}\to\rightarrow
H_k(\varGamma_1(N);{\frak k})\quad {\text{when $2 \leq k \not\equiv
0$}}\pmod{p-1} \\
e\,H_k(\varGamma_1(N)\cap \varGamma_0(p);{\frak r})/\varpi
\overset{\sim}\to\rightarrow
e\, H_k(\varGamma_1(N)\cap \varGamma_0(p);{\frak k}) \\
\qquad {\text{when $3  \leq k \not\equiv 0$}} \pmod{p-1}.
\endcases $$
\endproclaim
\demo{Proof}\!\! We clearly have a natural surjection from
$H_k(\varGamma_1(N);{\frak r})/\varpi$ to $H_k(\varGamma_1(N);{\frak k})$.
But the perfectness of the first pairing above, reduced modulo $\varpi$,
shows that the action of the former ring on $M_k(\varGamma_1(N);{\frak k})$
is faithful, which settles the first case. The same proof works in the
second case. \qed \enddemo

As before, we indicate by the superscript $``{}^{(i)}"$ the 
$\omega^i$-eigenspace with respect to the action of
$({\bold Z}/p{\bold Z})^{\times}$.
\proclaim{Corollary (1.4.3)} If $i \not\equiv -2 \pmod{p-1}$, the
pairing:
$$e\, {\Cal H}(N;{\frak r})^{(i)}
\times e\, M(N;\Lambda_{\frak r})^{(i)}
\to \Lambda_{\frak r}$$
defined by $(t, {\Cal F}):= {(\text{{\rm{the coefficient of $q$ in
${\Cal F}\mid t)$}}}}$ is perfect.
\endproclaim
\demo{Proof} Take a positive integer $d$ such that
$d \equiv i$ (mod $p-1$). Then the pairing in question reduced modulo
$\omega_d$ may be identified with the second one in (1.4.1) with
$k=d+2$.
\qed \enddemo

\demo{Remark {\rm{(1.4.4)}}} In [O4], 3.3, we discussed this type of
duality for Eisenstein components under the assumption $p\nmid
\varphi (N)$.
For such a component attached to
a pair $(\vartheta ,{\pmb{1}})$ with $\vartheta \not= \omega^{-2}$, we
remarked that its {\it a priori} proof simplifies the proof
of the main theorem of [O4]. The corollary above supplies such a proof for
$\vartheta$ satsisfying $\vartheta \mid_{({\bold Z}/p{\bold Z})^{\times}}
= \omega^i$ with $i\not\equiv -2 \pmod{p-1}$. 
In the next subsection, we will also give a direct proof for 
the remaining Eisenstein components when $p \nmid \varphi (N).$
\enddemo

\vskip3mm
\flushpar
{\bf 1.5. Eisenstein components.} Let $\vartheta$ and $\psi$ be 
primitive Dirichlet
characters of conductors $u$ and $v$, respectively.
We assume that $uv=N$ or $Np$, and that the condition (1.2.4) is satisfied
(so that we necessarily have $c=1$). We will write $E_k(\chi , \psi )$
for $E_k(\chi ,\psi ;1)$ in the following.  
 
We assume that ${\frak r}$ contains the values of
$\vartheta$ and $\psi$.
We set:
$$\vartheta = \chi \omega^i \tag 1.5.1$$
with a character $\chi$ whose conductor $u_0$ is prime to $p$. 
Let ${\frak M}={\frak M}(\vartheta ,\psi )$ be the Eisenstein
maximal ideal of $e\, {\Cal H}(N;{\frak r})$ associated 
with ${\Cal E}(\vartheta ,\psi )$ (1.2.9).
It is then easy to see that 
$e\, {\Cal H}(N;{\frak r})_{\frak M}$ is
contained in $e\, {\Cal H}(N;{\frak r})^{(i)}$. 

In the following discussions, we take and fix an integer
$d \geq 1$ congruent to $i$ modulo $p-1$,
and set $k:=d+2$.
It follows from the remark above that
(1.2.11) induces a ring isomorphism:
$$
e\,{\Cal H}(N;{\frak r})_{\frak M}/\omega_d \overset{\sim}\to\rightarrow
e\, H_k(\varGamma_1(N)\cap \varGamma_0(p);{\frak r})_{\frak M}.
\tag 1.5.2 $$
Here, we are identifying ${\frak M}$ with the corresponding maximal ideal
of $e\, {\Cal H}(N;{\frak r})^{(i)}$, and considering $e\, H_k(\varGamma_1(N)
\cap \varGamma_0(p);{\frak r})$ as a module over this ring.
By (1.2.7), the $({\pmb{1}},d)$-specialization of ${\Cal E}
(\vartheta ,\psi )$ is $E_k(\chi_1, \psi )$. On the other hand, 
$E_k(\chi , \psi )$ belongs to $M_k(\varGamma_1(N);{\frak r})$.
Let ${\frak m}={\frak m}(k;\chi ,\psi )$ 
be the Eisenstein maximal ideal of $H_k(\varGamma_1(N);
{\frak r})$ associated with this series. It is thus generated by
all 
$T(n)-\sum_{0<t\mid n}\chi (t)\psi (n/t)t^{k-1}$ and $\varpi$.
We now consider the duality of the same type as in 1.4 for the local
components attached to ${\frak m}$ and ${\frak M}$:
\proclaim{Proposition (1.5.3)} Suppose that the following conditions
are not simultaneously satisfied:

$$\cases
i \equiv -2 \,\,({\roman{mod}}\,\, p-1); \\
{\text{the orders of $\chi$ and $\psi$ are powers of $p$}};\\
(u_0,v)=1 ;\\
{\text{every prime factor of $v$ is
congruent to one modulo $p$}}.
\endcases $$
Then no non-zero constant belongs to 
$M_k(\varGamma_1(N);{\frak k})_{\frak m}$, and the pairings:
$$\cases
H_k(\varGamma_1(N);{\frak r})_{\frak m} \times
M_k(\varGamma_1(N);{\frak r})_{\frak m}\to {\frak r} \\
e\,{\Cal H}(N;{\frak r})_{\frak M} \times e\, M(N;
{\Lambda}_{\frak r})_{\frak M}\to \Lambda_{\frak r}
\endcases$$
defined as in $(1.4.1)$ and $(1.4.3)$, respectively, are perfect.
\endproclaim
\demo{Proof} 
By (1.3.3) and (1.3.4), the spaces
$e\, M_k(\varGamma_1(N)\cap \varGamma_0(p);{\frak k})_{\frak M}$
and $M_k(\varGamma_1(N);$

\flushpar
${\frak k})_{\frak m}$ coincide.
Therefore, in view of the proofs of (1.4.1) and (1.4.3),
it is enough to show that no non-zero constant belongs to 
the latter space.

When $i \equiv  -2  \pmod{p-1}$, i.e. $k \equiv 0 \pmod{p-1}$,
a non-zero constant $a\in {\frak k}$ 
certainly belongs to $M_k(\varGamma_1(N); {\frak k})$. 
It is invariant
under the diamond operators and satisfies
$$a \mid T(l)=
\cases (1+l^{k-1})a \quad{\text{if $l\nmid N$}} \\
a\quad {\text{if $l \mid N$}}\endcases $$
for all prime numbers $l$.
Thus if this element belongs to $M_k(\varGamma_1(N);{\frak k})_{\frak m}$,
it must be annihilated by ${\frak m}$.
Our assumption clearly rules this possibility out.
\qed \enddemo

We have especially shown that, when $p \nmid \varphi (N)$, 
the second pairing above is always perfect
(since the case where
$(\vartheta ,\psi )= (\omega^{-2},{\pmb 1})$ is excluded by the condition
(1.2.4)). By the same argument as in 1.4, we also obtain the following:

\proclaim{Corollary (1.5.4)} Under the same hypothesis as in $(1.5.3)$, 
we have canonical ring isomorphisms:
$$\cases
H_k(\varGamma_1(N);{\frak r})_{\frak m}/\varpi 
\overset\sim\to\rightarrow
H_k(\varGamma_1(N);{\frak k})_{\frak m}\\
e\, H_k(\varGamma_1(N)\cap \varGamma_0(p); {\frak r})_{\frak M}/\varpi
\overset\sim\to\rightarrow
e\, H_k(\varGamma_1(N)\cap \varGamma_0(p);{\frak k})_{\frak M}
\endcases $$
and the rings in the right hand side are isomorphic. 
\endproclaim

\vskip5mm
\flushpar
{\bf \S2. Ulmer's pairing}

\vskip3mm \flushpar
{\bf 2.1. Selmer group ${\roman{Sel}}(K,V)$.} The purpose of this
section is to construct a Hecke equivariant bilinear
pairing between the spaces of modular forms and cusp forms over
finite fields ((2.5.1), (2.5.5)).
Such a pairing is obtained by twisting the one constructed by
Ulmer [U], and so we begin by recalling his result.

As in Section 1, we fix a prime number $p \geq 5$, and a positive integer
$N$ prime to $p$. For the moment, until 2.5, we assume that
$N\geq 5$. In what follows, we use the same terminology and convention
for modular curves as in Gross' paper [Gr]. Let $Y_1(N)_{/{\bold F}_p}$
denote the fine moduli scheme classifying the pairs $(E, \alpha )$
consisting of an elliptic curve $E$ and a closed
immersion of group schemes $\alpha : {\pmb \mu}_N \hookrightarrow
E_N$, over ${\bold F}_p$-schemes. Its smooth compactification
$X_1(N)_{/{\bold F}_p}$ actually parametrizes the pairs $(E,\alpha )$
with $E$ a generalized elliptic curve
and $\alpha: {\pmb \mu}_N\hookrightarrow E_N$ a morphism
whose image meets every irreducible component of every geometric fibre
(cf. [Gr], Section 2).
The points of $X_1(N)_{/{\bold F}_p}-Y_1(N)_{/{\bold F}_p}$ 
are called cusps. On the other hand, the functor which assigns to
each ${\bold F}_p$-scheme $S$ the $S$-isomorphism classes of
the triples
$(E, \alpha ,\beta )$ consisting of:
$$\cases 
{\text{an elliptic curve $E$;}}\\
{\text{a closed immersion $\alpha :{\pmb {\mu}}_N\hookrightarrow E_N$
of group schemes;}}\\
{\text{a closed immersion $\beta : {\pmb{\mu}}_p \hookrightarrow
E_p$ of group schemes}}
\endcases \tag 2.1.1$$ 
over $S$,
is represented by an \'etale covering of 
$Y_1(N)_{/{\bold F}_p}-\{ {\roman{supersingular\,\, points}}\}$.
Its smooth compactification $I_1(N)$ is the Igusa curve of level
$Np$ over ${\bold F}_p$ (cf. [Gr], Section 5). $I_1(N)$ is geometrically
irreducible over ${\bold F}_p$.
A point of $I_1(N)$ is called ordinary, 
supersingular, or cuspidal according as its image to $X_1(N)_{/{\bold F}_p}$
corresponds to an ordinary elliptic curve, supersingular
elliptic curve, or a cusp, respectively.

Fix a finite extension ${\frak k}$ of ${\bold F}_p$, and set $K:=
{\frak k}(I_1(N))$, the function field of $I_1(N)$ over
${\frak k}$. Let ${\Cal E}$ be the base change to Spec$(K)$ of the
universal elliptic curve on 
(an open subscheme of) $I_1(N)$. Then if we denote by 
$V: {\Cal E}^{(p)}\to {\Cal E}$ the Verschiebung, the morphism $\beta$
for ${\Cal E}$ determines a canonical isomorphism:
${\bold Z}/p{\bold Z} \cong {\roman{Ker}}(V)$
by the Cartier-Nishi duality.
Let Sel$(K,V)$ be the Selmer group
as defined in [U], Section 1. It is a subgroup of $H^1(K, {\roman{Ker}}(V))$ 
satisfying certain local conditions;
but what we actually need is the following explicit description:
First we note that $H^1(K,{\bold Z}/p{\bold Z})\cong K/\wp (K)$ 
via the Artin-Schreier theory, where $\wp (x):=x^p-x$. 
Then, Sel$(K,V)$ can be identified with the subgroup 
of $K/\wp (K)$ represented by $f \in K$ satisfying the following conditions:
$$\cases
{\text{$f \in \wp (K_v)$ if $v$ is a cuspidal place of $K$;}}\\
{\text{$f \in R_v +\wp (K_v)$ if $v$ is an ordinary place of $K$;}}\\
{\text{$f = f_v + \wp (g_v)$ with $f_v , g_v \in K_v$  such that
$v(f_v)> -p$ }}\\
{\text{\quad if $v$ is a supersingular place of $K$}}
\endcases \tag 2.1.2$$
([U], Proposition 3.1).
Here, $R_v$ denotes the ring of integers of $K_v$.

There is a canonical isomorphism $\langle \enskip \rangle_p$
of $({\bold Z}/p{\bold Z})^{\times}$ onto Aut$(I_1(N)/
X_1(N)_{/{\bold F}_p})$ which, on ordinary points, is given by
$$\langle b \rangle_p: (E, \alpha ,\beta )\mapsto
(E, \alpha , b\beta ).
\tag 2.1.3$$
Consequently, 
the Galois group Gal$(K/{\frak k}(X_1(N)_{/{\bold F}_p}))$ is 
also canonically
isomorphic to $({\bold Z}/p{\bold Z})^{\times}$. As usual,
we indicate by the superscript $``{}^{(i)}"$ the eigenspace on which
$({\bold Z}/p{\bold Z})^{\times}$ acts via $\omega^i$.  We then
recall that Sel$(K,V)^{(k-1)}$ is identified with the subgroup of
$(K/\wp (K))^{(k)}$ represented by $f \in K^{(k)}$ satisfying the above
conditions.

\vskip3mm \flushpar
{\bf 2.2. Ulmer's mapping.} In what follows, we fix an integer
$k$ such that $2 \leq k \leq p-1$, and set $k':= p+1-k$
so that $2 \leq k' \leq p-1$. In this subsection, we recall the
definition of the (additive but not necessarily ${\frak k}$-linear) mapping
$$M_{k'}(\varGamma_1(N); {\frak k})\to {\text{Hom}}_{\frak k}
(S_{k'}(\varGamma_1(N);{\frak k}), {\frak k})
\tag 2.2.1 $$
constructed by Ulmer ([U], Theorem 7.8).

In general, if
$\phi : E\to S$ is a generalized elliptic curve, 
the invertible sheaf $\underline{\omega}_E$
on $S$ is defined as
the ${\Cal O}_S$-dual of the relative Lie algebra
$\underline{\roman{Lie}}(E^{reg})$, where $E^{reg}$ is the
smooth locus of $\phi$. 
Thus $\underline{\omega}_E$ is simply $\phi_{*}\Omega_{E/S}^1$ if
$E$ is a genuine elliptic curve over $S$.
A (geometrically defined) modular
form ${\pmb{f}}$ of weight $n$ on $\varGamma_1(N)$ defined over ${\frak k}$
in the sense of Katz is a rule which assigns to every pair $(E, \alpha )$
($E$ is a generalized elliptic curve over a ${\frak k}$-scheme and
$\alpha: {\pmb \mu}_N\hookrightarrow E_N$ is as in 2.1) a section
${\pmb f}(E,\alpha )$
of $\underline{\omega}^{\otimes n}_E$,
compatibly with base changes and isomorphisms
(cf. [Ka1], [Gr]). Let $\underline{\omega}$ be the invertible
sheaf on $X_1(N)_{/{\bold F}_p}$
corresponding to the universal generalized elliptic curve
over it. Then the space of modular forms above can be identified
with the space $H^0(X_1(N)_{/{\frak k}}, \underline{\omega}^{\otimes n})$,
where
$X_1(N)_{/\frak k}:=X_1(N)_{/{\bold F}_p}\otimes_{{\bold F}_p}{\frak k}$.
 
If ${\pmb f}$ is an element of this space,
we can evaluate it at the Tate curve ${\bold G}_m/q^{\bold Z}$ 
over ${\frak k}((q))$ together
with the natural $Id_N :{\pmb \mu}_N\hookrightarrow 
({\bold G}_m/q^{\bold Z})_N$, and we have:

$${\pmb f}({\bold G}_m/q^{\bold Z}, Id_N)=f \cdot (\frac{dt}t)^{\otimes n}
\quad
{\roman{with}}\enskip f \in {\frak k}[[q]] \tag 2.2.2$$
where $t$ is the standard parameter on ${\bold G}_m$.
The correspondence ${\pmb f}\mapsto f$
establishes an isomorphism: $H^0(X_1(N)_{/{\frak k}}, 
\underline{\omega}^{\otimes n})\cong M_n(\varGamma_1(N);{\frak k})$
for all $n \geq 2$.

We set
$$\cases X:=I_1(N)
\otimes_{{\bold F}_p}{\frak k};\\
Y:=X-\{ {\text{supersingular points}} \} =:{\text{Spec}}(R).
\endcases \tag 2.2.3$$
Let $\pi : I_1(N)\to X_1(N)_{/{\bold F}_p}$ be the natural morphism,
and denote by the same symbol $\underline{\omega}$
the inverse image of the previous $\underline{\omega}$ by $\pi$. 
Then there is a canonical section $a \in H^0(I_1(N), \underline{\omega})$
as described in [Gr], Proposition 5.2. Now let
$f \in M_n(\varGamma_1(N); {\frak k})$ correspond to
${\pmb f}\in H^0(X_1(N)_{/{\frak k}}, \underline{\omega}^{\otimes n})$ 
($n \geq 2$). Then the correspondence
$f \mapsto {\pmb f}/a^n$
gives an isomorphism:
$$\widetilde{M}:=\bigcup_{n=0}^{\infty}M_n(\varGamma_1(N); {\frak k})
\overset{\sim}\to\rightarrow R
\tag 2.2.4$$
(cf. [Gr], Proposition 5.5). Here, we are considering the left hand side
as a subset of ${\frak k}[[q]]$, and note that it coincides with
the union of $M_n(\varGamma_1(N);{\frak k})$ for all $n \geq 2$.
We will henceforth
identify these two rings, and consider each element
of $M_n(\varGamma_1(N); {\frak k}$) as giving an element of $R \subseteq
K$. Let Sel$(K,V)_j$ be the subgroup of Sel$(K,V)$
represented by forms of filtration $\leq j$.

We now recall Ulmer's construction. 
First, we define the mapping
$$M_{k'}(\varGamma_1(N); {\frak k}) \to {\text{Sel}}(K,V)^{(k-1)}
/{\text{Sel}}(K, V)_{p-1}^{(k-1)}
\tag 2.2.5 $$
by sending $f \in M_{k'}(\varGamma_1(N); {\frak k})$ to the class
of $h := \theta^{k-1}f$, where $\theta = q(d/dq)$ is the
$\theta$-operator of Serre and Katz. 
In [U], Theorem 7.8, (d), Ulmer defined this mapping on
$S_{k'}(\varGamma_1(N);{\frak k})$ for $1 \leq k \leq p$. 
Here, we are assuming that $k-1>0$, in order to assure that
$h$ above is  a cusp form (cf. Katz, [Ka2], II). 
Thus $h$ satisfies the first local condition in
(2.1.2), as well as others, and determines a class in Sel$(K,V)^{(k-1)}$. 

The kernel of this mapping consists of the {\it forms with
companions}; i.e. those $f \in M_{k'}(\varGamma_1(N); {\frak k})$ 
such that $\theta^{k}f = \theta g$
with some $g  \in M_k(\varGamma_1(N);{\frak k})$.
This follows from the same argument as that
given in [U],
using the isomorphism:
Sel$(K,V)_{p-1}^{(k-1)} $
$\cong M_k(\varGamma_1(N);{\frak k})^{\wp -{\text{cusp}}}$
([U], Theorem 7.8, (b)). But we remark that the definition
of the right hand side should be the forms 
(considered as elements of $R$) whose values at every
cuspidal place $v$ belong to $\wp ({\text{the residue field of $v$}})$,
rather than $\wp ({\frak k})$.

On the other hand, we have an injection
$${\text{Sel}}(K,V)^{(k-1)}/{\text{Sel}}(K,V)_{p-1}^{(k-1)}
\hookrightarrow H^1(X, {\Cal O}_X)^{(k)}
\tag 2.2.6$$
defined as follows:
Let $f \in K$ represent a class $\overline{f} \in {\text{Sel}}
(K,V)^{(k-1)}$. Then, for each place $v$ of $K$, there are
$f_v, g_v \in K_v$ such that $f = f_v +\wp (g_v)$ and
$v(f_v)>-p$ (resp. $v( f_v)\geq 0$)
when $v$ is supersingular (resp. otherwise). The mapping above 
sends $\overline{f}$ to the cohomology class of
$(g_v)_{v} \in {\bold A}_K/(K+ \prod_v R_v) \cong H^1(X,{\Cal O}_X)$,
${\bold A}_K$ being the ring of adeles of $K$.

Next, the Serre duality  gives an isomorphism
$$H^1(X,{\Cal O}_X)^{(k)} \cong {\text{Hom}}_{\frak k}
(H^0(X,\Omega_{X/{\frak k}}^1)^{(k'-2)}, {\frak k}).
\tag 2.2.7$$
Here, we let $b \in ({\bold Z}/p{\bold Z})^{\times}$ act on differentials
via the pull-back $\langle b \rangle_p^*$
by $\langle b \rangle_p$, and note that 
$k'-2 \equiv -k \pmod{p-1}$.

Finally, we have Serre's isomorphism
$$H^0(X, \Omega_{X/{\frak k}}^1)^{(k'-2)} \cong
S_{k'}(\varGamma_1(N);{\frak k})
\tag 2.2.8$$
(cf. [Gr], Proposition 5.7).

Combining (2.2.5)-(2.2.8), we obtain the mapping (2.2.1). Note that
(Sel$(K,V)$ is just an abelian group and) the composite of
(2.2.5) and (2.2.6) is not ${\frak k}$-linear unless ${\frak k}={\bold F}_p$. 
So, we make the following
\demo{Definition {\rm{(2.2.9)}}} We define the pairing, called 
{\it Ulmer's pairing},
$$( \, \,  , \, )_{k'}: M_{k'}(\varGamma_1(N);{\frak k}) \times
S_{k'}(\varGamma_1(N); {\frak k}) \to {\frak k}$$
as the the one obtained from (2.2.1) when ${\frak k}={\bold F}_p$,
and as its ${\frak k}$-bilinear extension in general.
\enddemo
Thus the left kernel of this pairing is the ${\frak k}$-linear
span of the left kernel of the pairing over ${\bold F}_p$; i.e.
it consists precisely of the forms in $M_{k'}(\varGamma_1(N);
{\frak k})$ having companions. 

\vskip3mm \flushpar
{\bf 2.3. Relation with Hecke operators.} We are now going to look 
at the compatibility of Hecke operators with respect to the mappings
(2.2.4)-(2.2.8).

For the moment, we fix a prime number $l$ different from $p$, and 
consider the Hecke operator $T(l)$.
We first recall related algebraic correspondences (cf. [Gr], Sections 3
and 5).
When $l$ does not divide $N$, 
let $Y(l)^0$ be the fine moduli scheme parametrizing quadruples 
$(E, \alpha ,
\beta , C)$, where $(E, \alpha , \beta )$ is as in (2.1.1)
over a ${\frak k}$-scheme, and
$C$ is a locally free subgroup scheme of $E_l$ of rank $l$. 
There is an \'etale morphism $Y(l)^0 \to Y-\{{\text {cusps}}\}$,
sending $(E, \alpha ,\beta , C)$ to $(E, \alpha , \beta )$,
which uniquely extends to a finite flat morphism
$\pi_1 : X(l) \to X$ with $X(l)$ normal.
On the other hand, 
we can associate
with $(E, \alpha , \beta ,C)$ 
a triple $(E' , \alpha^{\prime}, \beta^{\prime})$ by
$$\cases
E':= E/C; \\
\alpha^{\prime}:=\varphi \circ \alpha; \\
\beta^{\prime}:=\varphi \circ \beta
\endcases \tag 2.3.1 $$
with $\varphi : E\to E'$ the quotient morphism.
Let $\pi_2 : X(l) \to X$ be the morphism
which, on ordinary points, sends
$(E, \alpha ,\beta , C)$ to $(E', \alpha^{\prime},
\beta^{\prime})$. 
When $l$ divides $N$, we define $Y(l)^0$ as the fine moduli scheme
parametrizing quadruples as above 
such that $C \cap {\text{Im}}(\alpha )=0$, and then define
$\pi_1 :X(l) \to X$ and $\pi_2 : X(l)\to X$ in a similar manner.
Set 
$$Y(l):=X(l)-\{{\text{supersingular points}}\}. \tag 2.3.2$$

\proclaim{Proposition (2.3.3)} The notation being as above, the
following diagram commutes:
$$\CD
\widetilde{M} @>{\text{(2.2.4)}}>{\sim}> R \\
@V{lT(l)}VV @VV{\pi_{1 *}\circ \pi_2^*}V \\
\widetilde{M} @>{\sim}>{\text{(2.2.4)}}> R
\endCD $$
where $\pi_{1*}$ is the trace mapping from $H^0(Y(l),
{\Cal O}_{Y(l)})$ to $R$.
\endproclaim
\demo{Proof} We may identify an ordinary geometic point $P$ of $Y$
with a triple $(E, \alpha ,\beta )$ as in (2.1.1). Then for
$x \in R$, we have:
$$\align
\pi_{1*}\circ \pi_2^*(x)(P)&=\pi_{1*}\circ \pi_2^*(x)(E, \alpha ,\beta )
= \underset {Q=(E, \alpha ,\beta , C)\mapsto P}\to\sum \pi_2^*(x)(E, \alpha ,
\beta ,C) \\
&= \sum_{Q=(E, \alpha , \beta , C)\mapsto P}
x(E', \alpha^{\prime} ,\beta^{\prime}). \endalign$$
Here, the sums run over all the geometric points 
$Q=(E, \alpha ,\beta ,C)$ of $Y(l)$ 
such that $\pi_1(Q)=P$, and $(E',\alpha^{\prime},\beta^{\prime})$
is obtained from $Q$ by (2.3.1).

Take $f =\sum_{n=0}^{\infty}a_nq^n
\in M_k(\varGamma_1(N); {\frak k})$, and let ${\pmb{f}}$
be the geometrically defined modular form corresponding to $f$
by (2.2.2). 
The mapping (2.2.4) sends $f$ to $x:={\pmb{f}}/a^k$.
What we want to show is that
the value $\pi_{1*}\circ \pi_2^*(x)(P)$ for $P=
({\bold G}_m/q^{\bold Z}, Id_N, Id_p)$ 
coincides with $lf\mid T(l)$,
where ${\bold G}_m/q^{\bold Z}$ is the Tate curve over ${\frak k}((q))$
and $Id_p: {\pmb{\mu}}_p\hookrightarrow ({\bold G}_m/q^{\bold Z})_p$ 
is the natural embedding.

To do this, we may assume that ${\frak k}$ contains a primitive
$l$-th root of unity $\zeta_l$. 
We first treat the case where $l \nmid N$.
Then there are $l+1$ subgroups
$$\cases {\text{(i)}}\,\, C_0:={\pmb{\mu}}_l=\langle \zeta_l \rangle ;\\
{\text{(ii)}}\,\,C_i:=\langle \zeta_l^i q^{1/l}\rangle 
\,\,(i=1,\cdots,l)
\endcases $$
of order $l$
in ${\bold G}_m/q^{\bold Z}$ over ${\frak k}((q))$. Let 
$(E_i^{\prime}, \alpha_i^{\prime}, \beta_i^{\prime})$ be the
object corresponding to
$({\bold G}_m/q^{\bold Z},
Id_N, Id_p, C_i)$ by (2.3.1).
Noting that $a^k({\bold G}_m/q^{\bold Z}, Id_N,
Id_p)=(dt/t)^{\otimes k}$ and $a^k({\bold G}_m/q^{\bold Z}, lId_N,
lId_p)=l^{-k}(dt/t)^{\otimes k}$ (cf. [Gr], Proposition 5.2),
a direct calculation as in Katz [Ka1], 1.11 shows that:
$$\cases
{\text{(i)\,\, $x(E_0^{\prime}, \alpha_0^{\prime} ,\beta_0^{\prime})
= l^k \sum_{n=0}^{\infty}b_nq^{nl}$
if $f \mid \langle l \rangle =\sum_{n=0}^{\infty}
b_nq^n$;}} \\
{\text{(ii)}}\,\, x(E_i^{\prime},\alpha_i^{\prime} ,\beta_i^{\prime})
= \sum_{n=0}^{\infty} a_n(\zeta_l^iq^{1/l})^n
\,\,(i=1,\cdots ,l).
\endcases $$
The sum of the
values in the case (ii) is equal to $l\sum_{n=0}^{\infty}a_{nl}q^n$,
which completes the proof when $l\nmid N$.

When $l\mid N$, the subgroups to be considered are those in the 
case (ii), and the computation above shows that our claim also
holds in this case.
\qed
\enddemo

It is easy to see that $\pi_{1*}\circ \pi_2^*$ preserves
$\wp (R)$ and Sel$(K, V)^{(k-1)} \subseteq R/\wp (R)$. The proposition 
above and the well-known property:
$$T(n)\circ \theta = n\theta \circ T(n) \tag 2.3.4$$
show that the left square in the following diagram commutes:
$$\CD
M_{k'}(\varGamma_1(N);{\frak k}) @>(2.2.5)>>
{\text{Sel}}(K,V)^{(k-1)}/{\text{Sel}}(K,V)_{p-1}^{(k-1)}
@>(2.2.6)>> H^1(X,{\Cal O}_X)^{(k)}\\
@Vl^kT(l)VV @VV\pi_{1*}\circ \pi_2^*V
@VV\pi_{1*}\circ \pi_2^*V \\
M_{k'}(\varGamma_1(N);{\frak k}) @>>(2.2.5)> 
{\text{Sel}}(K,V)^{(k-1)}/{\text{Sel}}(K,V)_{p-1}^{(k-1)}
@>>(2.2.6)> H^1(X,{\Cal O}_X)^{(k)}.
\endCD \tag 2.3.5 $$
Also from the definition of the
mapping (2.2.6) recalled in 2.2, the right square 
commutes.
Here, $\pi_{1*}$ in the right vertical arrow
means the trace mapping for the cohomology.
Since the canonical mapping and the trace mapping for the
cohomology groups
are transformed to each other by the Serre duality, we 
have the following commutative diagram:
$$\CD
H^1(X,{\Cal O}_X)^{(k)} @>(2.2.7)>>{\text{Hom}}_{\frak k}
(H^0(X,\Omega_{X/{\frak k}}^1)^{(k'-2)} ,{\frak k})\\
@V\pi_{1*} \circ \pi_2^* VV @VV{}^t(\pi_{2*} \circ \pi_1^*)V \\
H^1(X,{\Cal O}_X)^{(k)} @>>(2.2.7)>{\text{Hom}}_{\frak k}
(H^0(X,\Omega_{X/{\frak k}}^1)^{(k'-2)}, {\frak k})
\endCD \tag 2.3.6$$
where the superscript $``{}^t"$ means the transpose. 

Finally, as shown by Gross ([Gr], Poposition 5.9), via the
isomorphism (2.2.8), $T(l)$ on the right commutes with $\pi_{1*}
\circ \pi_2^*$ on the left. 

By a similar argument, we can treat the diamond
operators: For $a \in ({\bold Z}/N{\bold Z})^{\times}$, let
$\langle a \rangle_N $ be
the automorphism of $X$ which is given by
$$\langle a \rangle_N : (E,\alpha ,\beta )\mapsto
(E, a\alpha ,\beta )
\tag 2.3.7$$
on ordinary points.
Then a similar but simpler argument as (2.3.3) shows that,
via (2.2.4), the diamond operator $\langle a \rangle $ on 
$\widetilde{M}$ corresponds to the ring automorphism of $R$
induced by $\langle a \rangle_N$. Proceeding
as above, we conclude that via (2.2.5)-(2.2.7),
$\langle a \rangle $ on $M_{k'}(\varGamma_1(N);{\frak k})$
corresponds to
${}^t\langle a^{-1} \rangle_N^*$ on ${\text{Hom}}_{\frak k}
(H^0(X,\Omega_{X/{\frak k}}^1)^{(k'-2)},{\frak k})$, the
transpose of the mapping taking the pull-back of differentials. 
So, again by [Gr], Proposition 5.9, 
$\langle a \rangle$ on the left
and ${}^t\langle a^{-1} \rangle $ on the right commute through
(2.2.1).

We have especially seen that the pairing (2.2.9) 
is {\it not} compatible with Hecke operators.

\vskip3mm \flushpar
{\bf 2.4. Operator $w_N$.} To remedy the situation in the previous
subsection, we will $``$twist" Ulmer's pairing (2.2.9) by the operator
$w_N$ to obtain a
Hecke equivariant one, in quite an analogous way as one obtains the
$``$twisted Weil pairing" from the usual Weil pairing.

To do this, we assume in this subsection that ${\frak k}$
contains a primitive $N$-th root of unity $\zeta_N$; 
and we throughout fix this choice of $\zeta_N$
(on which our operator
$w_N$ depends). For an elliptic curve $E$ over a ${\frak k}$-scheme,
we denote by
$$e_{N,E}(\,\, ,\,):E_N \times E_N \to {\pmb{\mu}}_N
\tag 2.4.1 $$
the Weil pairing.
Let $w_N$ be the automorphism of $X$ 
which sends an ordinary triple $(E,\alpha ,\beta )$ to $(E^*,
\alpha^*, \beta^*)$ defined by:
$$\cases
{\text{$E^*:=E/\alpha ({\pmb{\mu}_N})$, with the quotient morphism
$\phi :E \to E^*$;}}\\ 
{\text{$\alpha^*(\zeta_N):=\phi (t_{\alpha})$ with a section $t_{\alpha}$
of $E_N$ such that $e_{N,E}(\alpha (\zeta_N), t_{\alpha})=\zeta_N$;}}\\
\beta^*:=\phi \circ \beta .
\endcases \tag 2.4.2 $$
One easily checks that
$$w_N^2=\langle -1 \rangle_N \circ
\langle N \rangle_p 
\tag 2.4.3$$
in the same manner as [Gr], Proposition 6.7, 3).

\proclaim{Proposition (2.4.4)} Let $l$ be a prime number different
from $p$, and let $\pi_1,\pi_2 : X(l) \to X$ be as in the previous
subsection.
Then as algebraic correspondences on
$X\times X$, we have:
$$w_N^{-1}\circ \pi_{2*}\circ \pi_1^* \circ w_N =
\langle l \rangle_p  \circ \pi_{1*}\circ \pi_2^*.$$
\endproclaim
\demo{Proof} We give the proof only in the case where
$l \mid N$, since the proof when $l \nmid N$ is similar but simpler.

We first consider the image of an ordinary point
$P=(E, \alpha ,\beta )$ under the correspondence $w_N \circ \pi_{2*}
\circ \pi_1^* \circ w_N$.
Let $w_N(P):=P^*=(E^* ,\alpha^*, \beta^*)$ and $\phi$ be as above. 
If we denote by $C_i$ ($ 1\leq i \leq l$) the finite
subgroup schemes
of $E^*$ of order $l$ not contained in the image of $\alpha^*$, the
image of $P^*$ under $\pi_{2*}\circ \pi_1^*$ is the sum of
$P_i:=(E_i ,\alpha_i,\beta_i)$ with $E_i:=E^*/C_i$,
$\varphi_i:E^* \to E_i$ the quotient morphism,
$\alpha_i:=\varphi_i\circ \alpha^*$ and $\beta_i:=\varphi_i \circ
\beta^*$.
Set $w_N(P_i) =(E_i^*, \alpha_i^*,\beta_i^*)$,
and let $\pi_i: E_i \to E_i^*=E_i/\alpha_i(\pmb{\mu}_N)$
be the quotient morphism. Since
$\alpha_i({\pmb{\mu}}_N)
=\varphi_i \circ \phi (E_N)$, there is a unique isogeny $\varphi_i^*:
E \to E_i^*$ which makes the following diagram commutative:
$$\CD
E @>\phi>>E^* @>{}^t\phi>> E \\
@. @V\varphi_iVV @VV\varphi_i^*V \\
@. E_i @>>\pi_i> E_i^*.
\endCD $$
Here, ${}^t\phi$ is the isogeny dual to $\phi$.
Since $C_i\cap \alpha^*({\pmb{\mu}}_l)=0$, we have ${}^t\phi
(C_i)=
{\text{Ker}}(\varphi_i^*).$
Also, since $\phi \circ {}^t\phi (C_i)=0$,
we have ${}^t\phi (C_i)=\alpha ({\pmb{\mu}}_l)$,
which is independent of $i$. 
Therefore, if we denote by $\psi : E \to E/\alpha ({\pmb{\mu}}_l)$ the
quotient morphism, there is a unique isomorphism $\eta_i : E_i^*
\to E/\alpha ({\pmb{\mu}}_l)$ satisfying $\psi = \eta_i \circ \varphi_i^*$.
On the other hand, for a point $t$ of $E_{iN}$, we have:
$$e_{N,E_i}(\alpha_i(\zeta_N), t) =
e_{N,E^*}(\alpha^*(\zeta_N),{}^t\varphi_i(t))$$
$$=e_{N,E^*}(\phi (t_{\alpha}),{}^t\varphi_i(t))
=e_{N,E}({}^t\phi \circ {}^t\varphi_i (t), t_{\alpha})^{-1}$$
where $t_{\alpha}$ is a point of
$E_N$ as in (2.4.2). Noting that ${}^t\phi (E_N^*)
=\alpha ({\pmb{\mu}}_N)$, we see that
a point 
$t_{\alpha_i}$ 
of $E_{iN}$ satisfies
$e_{N,E_i}(\alpha_i(\zeta_N),t_{\alpha_i})=\zeta_N$ if
and only if ${}^t\phi \circ {}^t\varphi_i(t_{{\alpha}_i})=
\alpha (\zeta_N^{-1})$.

By (2.4.3), we conclude that the image of
$P=(E, \alpha ,\beta )$ under $w_N^{-1}\circ \pi_{2*}\circ \pi_1^* \circ
w_N$ is the sum of $(E':=E/\alpha ({\pmb{\mu}}_l), \alpha_i^{\prime},
\psi \circ \beta )$ ($1 \leq i \leq l$) with $\alpha_i^{\prime}$ defined as
follows: Let $t_{\alpha_i}^{\prime}$
be a point of $E_{iN}$ 
such that ${}^t\phi \circ {}^t\varphi_i(t_{\alpha_i}^{\prime})=
\alpha (\zeta_N)$. Then $\alpha_i^{\prime}$ sends $\zeta_N$ to
$\eta_i \circ \pi_i(t_{\alpha_i}^{\prime})$.

We next consider $\pi_{1*}\circ \pi_2^*(P)$.
$\pi_2^*(P)$ is a
sum of quadruples $(E'', \alpha'' ,\beta'', C)$ satisfying the following
conditions: Let $\xi : E'' \to E''/C$ be the quotient morphism. Then there
is an isomorphism of
$E''/C$ to $E$ through which $\xi \circ \alpha^{\prime \prime}$ 
(resp. $\xi \circ \beta^{\prime \prime}$)
corresponds to $\alpha$ (resp. $\beta$). We may
identify $E''/C$ with $E$,
Ker$({}^t\xi )$ with $\alpha ({\pmb{\mu}}_l)$,
and hence $\xi :E''\to E''/C$
with ${}^t\psi : E' \to E$, using the notation above.
$\pi_{1*}\circ \pi_2^*(P )$ is thus a sum of $(E', \alpha^{\prime},
\beta^{\prime})$
satisfying: (i) ${}^t\psi \circ \alpha^{\prime}=\alpha$ and 
(ii) ${}^t\psi \circ
\beta^{\prime} = \beta$.
The second condition
implies that $\beta^{\prime}=l^{-1}\psi \circ \beta$. 
As for $\alpha^{\prime}$,
if we fix one $\alpha^{\prime}$ satisfying (i), all $\alpha_t^{\prime}$ sending
$\zeta_N$ to $\alpha^{\prime} (\zeta_N) +t$ ($t \in {\text{Ker}}({}^t\psi )$)
also satisfy (i). Noting that $(E', \beta^{\prime})$ 
has no non-trivial automorphism,
we conclude that $\pi_{1*}\circ \pi_2^*(P)$ is the sum of
(non-isomorphic) $l$ triples $(E', \alpha_t^{\prime}, l^{-1}\psi \circ
\beta )$ for $t \in {\text{Ker}}({}^t\psi )$.

We now return to the situation of the first part. From the relation:
$\psi \circ {}^t\phi = \eta_i \circ \pi_i\circ \varphi_i$, we have
${}^t\phi \circ {}^t\varphi_i= {}^t\psi \circ \eta_i \circ \pi_i$.
It follows that ${}^t\psi \circ \alpha_i^{\prime}(\zeta_N)=
{}^t\phi \circ {}^t\varphi_i(t_{\alpha_i}^{\prime})=\alpha (\zeta_N)$,
so that $\alpha_i^{\prime}$ satisfies the condition (i) above.
Since the image of a point under $w_N^{-1} \circ \pi_{2*} \circ 
\pi_1^* \circ w_N$ consists of $l$ non-isomorphic triples 
generically (because the same holds for $\pi_{2*}\circ \pi_1^*$),
our conclusion follows.
\qed
\enddemo
\proclaim{Corollary (2.4.5)} For any differential form $\omega$
on $X$, we have:
$$w_N^*\circ \pi_{1*}\circ \pi_2^* \circ w_N^{-1*}(\omega )
=\pi_{2*}\circ \pi_1^*\circ \langle l \rangle_p^*(\omega ).$$
\endproclaim

By a similar but much simpler argument as (2.4.4), one obtains
$${\text{$w_N^{-1}\circ \langle a \rangle_N \circ w_N
= \langle a^{-1} \rangle_N$  for all $a \in ({\bold Z}/N{\bold Z})^{\times}$}}
\tag 2.4.6 $$
so that for any differential form $\omega$ on $X$, we have:

$$w_N^*\circ \langle a \rangle_N^*\circ w_N^{-1*}(\omega )
=\langle a^{-1} \rangle_N^*(\omega ).
\tag 2.4.7 $$
It is also easy to see that $w_N$ commutes with $\langle b \rangle_p$.

\vskip3mm \flushpar
{\bf 2.5. Twisting Ulmer's pairing.} So far, we argued under the assumption
that $N \geq 5$. We now drop this assumption in the following
\proclaim{Theorem (2.5.1)} Let $k'$ be an integer such that $2 \leq k'
\leq p-2$, and $N$ an arbitrary positive integer. Assume that ${\frak k}$
contains a primitive $N$-th root of unity.
Then there is a ${\frak k}$-bilinear pairing
$$(\, \, , \,)_{k'}^*: M_{k'}(\varGamma_1(N);{\frak k}) \times
S_{k'}(\varGamma_1(N);{\frak k}) \to {\frak k}$$
with respect to which all Hecke operators and diamond operators
are self-adjoint. Moreover,
the left kernel of this pairing consists of the forms with
companions.
\endproclaim

To prove this theorem, we assume for the moment that $N\geq 5$.
By the final remark in 2.4, $w_N^*$ induces a
${\frak k}$-linear automorphism of
$H^0(X,\Omega_{X/{\frak k}}^1)^{(k'-2)}$. 
Let us denote by $``\mid w_N"$ the automorphism of
$S_{k'}(\varGamma_1(N); {\frak k})$ corresponding to $w_N^*$
via (2.2.8).
We then set
$$(f , g )_{k'}^*:=(f,g \mid w_N)_{k'} \tag 2.5.2$$
where the right hand side is Ulmer's pairing (2.2.9). 
Since we just altered the right variable by applying
a ${\frak k}$-linear automorphism, this new pairing has the 
same left kernel as Ulmer's one. Hence the left kernel
consists of the forms with
companions, as remarked at the end of 2.2.
Also, it is 
clear from the argument in 2.3 and 2.4 that
$$(f\mid T,g )_{k'}^*=(f, g \mid T)_{k'}^* \tag 2.5.3$$
holds for $T=T(l)$ with a prime number $l \not= p$, or for
$T=\langle a \rangle$. We note that the argument up to this point
applies equally 
to the case where $k'=p-1$.
The proof of
(2.5.1) when $N \geq 5$ will be complete by the following

\proclaim{Lemma (2.5.4)} Let $k'$ satisfy $2 \leq k' \leq p-2$.
For any $N \geq 1$ and any finite field ${\frak k}$ of characteristic $p$,
the Hecke algebra $H_{k'}(\varGamma_1(N);{\frak k})$ is
generated over ${\frak k}$ by all $T(l)$ with 
prime numbers $l \not= p$ and the
diamond operators; or equivalently
by all $T(n)$ with $n$ prime to $p$.
\endproclaim

\demo{Proof} It is obvious that the two ${\frak k}$-subalgebras
generated by the indicated operators are the same.
Call this subalgebra $H'$.
By the duality (1.4.1) and (1.4.2), we know that 
$M_{k'}(\varGamma_1(N);{\frak k})$ is ${\frak k}$-dual to
$H_{k'}(\varGamma_1(N);{\frak k})$, and our assertion is
equivalent to the injectivity of the mapping:
$$M_{k'}(\varGamma_1(N);{\frak k}) \to {\roman{Hom}}_{\frak k}
(H', {\frak k})$$
which sends $f$ to the mapping: $t \mapsto a(1;f\mid t)$.
If $f$ belongs to the kernel of this mapping,
it must be a power series in $q^p$. The injectivity of 
the $\theta$-operator (cf. [Ka2], II) then
implies that $f=0$. (We have followed the proof of Ribet's lemma
given in Wiles [Wi], Chapter 2, \S2.)
\qed \enddemo

\demo{Variant {\rm{(2.5.5)}}} We consider here the case where $k'=p-1$,
other hypotheses being as above.
Set $H:=H_{p-1}(\varGamma_1
(N);{\frak k})$ and let $H'$ be its ${\frak k}$-subalgebra generated by
the operators indicated in (2.5.4).
 
Take Dirichlet characters $\chi$ and $\psi$ 
satisfying (1.2.2) with $u_0v=N$, and assume that
$\chi$ and $\psi$ do not satisfy the conditions in (1.5.3) simultaneously.
Let $\overline{\frak m}$ (resp. $\overline{\frak n}$) be the 
Eisenstein maximal ideal of $H$ (resp. $H'$) corresponding to
$E_{p-1}(\chi ,\psi )$ (cf. 1.5), so that
$\overline{\frak n}=H' \cap \overline{\frak m}$. 
The proof of (1.5.3) assures us that
$M_{p-1}(\varGamma_1(N);{\frak k})_{\overline{\frak n}} \cap {\frak k}=
\{ 0\}$, and hence we obtain a perfect pairing:
$$H_{\overline{\frak n}}\times 
M_{p-1}(\varGamma_1(N);{\frak k})_{\overline{\frak n}}
\to {\frak k}$$
where we are considering the localizations of $H'$-modules.
The proof of (2.5.4) then shows that $H_{\overline{\frak n}}^{\prime}
=H_{\overline{\frak n}}$, since the kernel of $\theta$ on
$M_{p-1}(\varGamma_1(N);
{\frak k})$ consists of constants.
We conclude that the 
restriction of the pairing (2.5.2)
$$(\, \, ,\, )_{p-1}^*:
M_{p-1}(\varGamma_1(N);{\frak k})_{\overline{\frak m}}\times
S_{p-1}(\varGamma_1(N);{\frak k})_{\overline{\frak m}}\to
{\frak k}$$
enjoys all the properties required in (2.5.1).
\enddemo

We now turn to the case where $N \leq 4$.
As remarked in [U], pages 253 and 263, this case can be
handled by adding an auxiliary level structure and taking
invariants. We explain this in some details below.
Take an integer $m\geq 3$ prime to $Np$.
We consider the moduli problem classifying the quadruples $(E, \alpha ,
\beta ,\gamma )$ over ${\bold F}_p$-schemes, where $E$, $\alpha$ and $\beta$
are the same as in (2.1.1), and $\gamma$ is a full level $m$ structure
$$\gamma : {\bold Z}/m{\bold Z}\times {\bold Z}/m{\bold Z}
\overset{\sim}\to\rightarrow E_m.$$
It is represented by a smooth and affine 
curve ${\Cal Y}_m$ over ${\bold F}_p$.

Let $\Phi_m(X)$ be the reduction modulo $p$
of the $m$-th cyclotomic polynomial, and set 
${\pmb{\mu}}_m^{\times}:={\roman{Spec}}\, ({\bold F}_p[X]/
(\Phi_m(X)))$.
Then $\gamma$ and the Weil pairing induce
a morphism: ${\Cal Y}_m\to {\pmb{\mu}}_m^{\times}$
(cf. Katz and Mazur [KM], Chapter 9).
We fix an irreducible factor $P(X)\in {\bold F}_p[X]$ of $\Phi_m(X)$
and its root $\zeta_m\in {\overline{\bold F}}_p$, and consider ${\frak k}
:={\bold F}_p(\zeta_m)$ as the coordinate ring of an irreducible component of
${\pmb{\mu}}_m^{\times}$ via $X \mapsto \zeta_m$. Let $Y_m=
{\roman{Spec}}\, (R_m)$ be the inverse image of ${\Cal Y}_m$ to
Spec$\, ({\frak k})$. As an ${\frak k}$-scheme, $Y_m$ is geometrically
irreducible, and represents the functor
classifying the quadruples $(E,\alpha ,\beta ,\gamma )$
where $\gamma$ has determinant $\zeta_m$ (cf. [KM], loc. cit.).
All the
arguments in 2.2 applies to this situation, replacing $\varGamma_1(N)$
by $\varGamma_1(N)\cap \varGamma (m)$,
and we write (2.2.$*$)${}_m$
for the statement corresponding to (2.2.$*$) in this situation. 

Let $G$ be the subgroup of $GL_2({\bold Z}/m{\bold Z})$ consisting
of matrices whose determinants are powers of $p$.
As an ${\bold F}_p$-scheme, $Y_m$ classifies the
quadruples $(E, \alpha ,\beta ,\gamma )$ where the determinant of
$\gamma$ is a root of $P(X)$. 
So, $G$ acts on $Y_m$ 
(and hence on $R_m$) in such a
way that $g \in G$ sends $(E,\alpha ,
\beta ,\gamma )$ to $(E, \alpha ,\beta ,\gamma \circ g)$.
On the other hand, identifying
$M_{k'}(\varGamma_1(N)\cap \varGamma (m);{\frak k})$ 
with the space of geometrically defined modular forms,
we can let $g \in  G$ act on this 
${\frak k}$-vector space $\det g$-linearly by the rule:
$$g {\pmb{f}}(E, \alpha ,\gamma ):= {\pmb{f}}
(E, \alpha ,\gamma \circ g).$$
Then since the canonical section $a \in H^0(Y_m, \underline{\omega})$
is invariant under $G$,
the isomorphism (2.2.4)${}_m$ is $G$-equivariant.
Since $\theta$ commutes with the action of $G$
(which follows from the description of $\theta$ given in
[Gr], Proposition 5.8), (2.2.5)${}_m$ is compatible with the 
$G$-action. We conclude that the pairing 
$$(\, \, ,\, )_{k',m}:M_{k'}(\varGamma_1(N)\cap \varGamma (m)
;{\frak k})\times
S_{k'}(\varGamma_1(N)\cap \varGamma (m)
;{\frak k})\to {\frak k}$$
obtained by composing (2.2.5)${}_m$-(2.2.8)${}_m$ satisfies
$$(gf ,gh)_{k',m}=(f, h)_{k',m}^{\det g}
\quad {\text{for all $g$}}\in G. \tag 2.5.6$$
Thus taking invariants under $G$, we obtain
a pairing as (2.2.9) over ${\bold F}_p$ for $N \leq 4$.
It is easy to see that this pairing is independent of
the choice of $m$ up to multiplication by an element of ${\bold F}_p^{\times}$
as long as the order of $G$ is prime to $p$, i.e. no prime factor
of $m$ is congruent to $\pm 1$ (mod $p$). We fix such a choice of
$m$, and call this pairing $(\, \, ,\,)_{k'}$.

Let us see that the left kernel of this pairing consists of the
forms with companions when $N\leq 4$ also. 
First assume that $p \equiv 1$ (mod 4). When
$N=1$, we take $m=4$. If $f \in M_{k'}(\varGamma_1(1);
{\bold F}_p)$ belongs to the left kernel of $(\, \, ,\, )_{k'}$, 
it follows from (2.5.6)
above that $f$ has a companion $h \in M_k(\varGamma (4); {\bold F}_p)$.
From the commutativity of $\theta$ and the action of $G=SL_2({\bold Z}
/4{\bold Z})$, we
conclude that $h$ actually belongs to $M_k(\varGamma_1(1);
{\bold F}_p)$, which proves our claim for $N=1$. 
On the other hand, if we start with $N=1$, $m=4$ and take 
invariants under the image of $\varGamma_1(2)$ (resp.
$\varGamma_1(4)$) in $G$, we obtain from $(\, \, ,\, )_{k',4}$ the
pairing $(\, \, ,\, )_{k'}$ for $N=2$ (resp. $N=4$), once we fix an isomorphism
${\bold Z}/4{\bold Z}\cong {\pmb \mu}_4$ over ${\bold F}_p$.
Thus our claim for $N=2$ (resp. $N=4$) 
follows in the same way. 
Starting with $N=3$ and $m=4$, we also settle the case 
where $N=3$. When $p \not\equiv 1$ (mod 4), $p-1$ has an odd
prime factor $l$. So, starting with $N\leq 4$ and $m=l$, we can
argue in the same manner as above to conclude.

Finally, take an $m\geq 5$ whose prime factors are not congruent to
$\pm 1$ (mod $p$), and consider the Igusa curve $X^{(m)}$ of level $Nmp$ over
${\bold F}_p(\zeta_N)$. Its ordinary points parametrize
quadruples $(E,\alpha ,\beta ,\delta )$
with $\delta : {\pmb \mu}_m
\hookrightarrow E_m$. We can define an automorphism $w_N^{(m)}$
of $X^{(m)}$ over ${\bold F}_p(\zeta_N)$ by sending $(E,\alpha ,\beta ,
\delta )$ to $(E^*,\alpha^* ,\beta^*, \delta^*)$ with $\delta^*=
\phi \circ \delta$. Clearly, this automorphism commutes with 
$w_N$ defined as in 2.4 on $X$, the Igusa curve of level $Np$
over ${\bold F}_p(\zeta_N)$, via
the natural projection. The argument above shows that Ulmer's pairing
of level $Nm$ restricts to that of level $N$. Therefore,
it follows from the argument in 2.4 that the twisted pairing 
defined by (2.5.2) is compatible with $T(l)$ ($l \nmid mp$) and
the diamond operators. Since there are infinitely many prime
numbers $m \not\equiv \pm 1$ (mod $p$), the proof of
(2.5.1) and (2.5.5) is now complete.
 
\vskip5mm
\flushpar
{\bf \S3. Structure of $p$-adic Hecke algebras
of Eisenstein type}

\vskip3mm
\flushpar
{\bf 3.1. Companions in Eisenstein components.} In this section,
we fix a positive integer $d$
such that $1 \leq d \leq p-3$, and set
$k:=d+2$ and $k':=p+1-k$. 

Let $\chi$ and $\psi$
be primitive Dirichlet characters of conductors $u_0$ and $v$, 
respectively. We assume that $(\chi \psi )(-1)=(-1)^k(=
(-1)^{k'})$ and $u_0v=N$. We also assume that the conditions in
(1.5.3) are not simultaneously satisfied with $i\equiv d$
(mod $p-1$). (We will
soon assume that $p$ does not divide $\varphi (N)$ so that this
assumption only excludes the case where $k=p-1$ and
$\chi = \psi ={\pmb{1}}$.)

As before, we let ${\frak r}$ be the ring of integers of
a finite extension of ${\bold Q}_p$ containing the values of
$\chi$ and $\psi$, and ${\frak k}={\frak r}/(\varpi )$ 
its residue field.

Let $E_k(\chi ,\psi )= E_k(\chi ,\psi ;1)
\in M_k(\varGamma_1(N);{\frak r})$ be the 
Eisenstein series defined by (1.2.1). By abuse of notation, we use
the same symbol to denote its reduction modulo $\varpi$ lying in
$M_k(\varGamma_1(N);{\frak k})$.

\proclaim{Lemma (3.1.1)} As modular forms {\rm{(mod}} $p)$, we
have:
$$\theta^{k'}E_k(\chi ,\psi )=\theta E_{k'}(\psi ,\chi )$$
i.e. $E_k(\chi ,\psi )$ and $E_{k'}(\psi ,\chi )$ are companions
to each other.
\endproclaim
\demo{Proof} Straightforward. \qed \enddemo

As in 1.5, we let ${\frak m}:={\frak m}(k; \chi ,\psi) \subset
H_k(\varGamma_1(N);{\frak r})$ be the maximal ideal associated with
$E_k(\chi ,\psi )$. Also, we let ${\frak m}^{\prime} 
:={\frak m}(k';\psi ,\chi ) \subset
H_{k'}(\varGamma_1(N);{\frak r})$ be the maximal ideal corresponding
to $E_{k'}(\psi ,\chi )$.

\proclaim{Proposition (3.1.2)} Let the notation and the hypotheses 
be as above, and
suppose that $f \in M_k(\varGamma_1(N);
{\frak k})$ has a companion $g \in M_{k'}(\varGamma_1(N);
{\frak k})$: $\theta^{k'}f=\theta g$.
If $f$ belongs to $M_k(\varGamma_1(N);{\frak k})_{\frak m}$,
then $g$ belongs to $M_{k'}(\varGamma_1(N);
{\frak k})_{{\frak m}^{\prime}}$. The converse also holds when
$k\not= p-1$.
\endproclaim
\demo{Proof} Since $2 \leq k' \leq p-2$,
$H_{k'}(\varGamma_1(N);{\frak k})$ is generated by
all $T(n)$ with $n$ prime to $p$ by (2.5.4), and hence the image 
of ${\frak m}'$ 
in $H_{k'}(\varGamma_1(N);{\frak k})$ is generated by all
$$\eta^{\prime} (n):= T(n)-\sum_{0<t\mid n}\psi (t)\chi (\frac{n}{t})t^{k'-1}$$
with $n$ prime to $p$.

Set
$$\eta (n):=T(n)-\sum_{0<t \mid n}\chi (t)\psi (
\frac{n}{t})t^{k-1}$$
for $n \not\equiv 0$ (mod $p$), which belongs to the image of
${\frak m}$ in $H_{k}(\varGamma_1(N);{\frak k})$.

From the
relation (2.3.4), one easily checks that
$$\theta^{k'-1}\circ \eta (n)^m=(n^{-k'+1}\eta^{\prime}(n))^m
\circ \theta^{k'-1}.$$ Thus if $f$ belongs to $M_k(\varGamma_1(N);
{\frak k})_{\frak m}$, we have:
$$0=\eta^{\prime}(n)^m\circ \theta^{k'-1}f=\eta^{\prime}(n)^m\circ
\theta^{p-1}g= \theta^{p-1}\circ \eta^{\prime}(n)^mg$$
for $m$ sufficiently large. This shows that $g$ is annihilated by
a power of
${\frak m}^{\prime}$, and hence settles the first assertion.

When $k\not= p-1$, we can reverse the roles of $f$ and $g$, since
$\theta^{k'}f=\theta g$ is equivalent to $\theta^kg=\theta f$.
\qed \enddemo

In the situation above, the companion $g$ of $f$ is unique; and
our hypothesis on $\chi$ and $\psi$ guarantees that the
correspondence $f\mapsto g$ is injective on
$M_k(\varGamma_1(N);{\frak k})_{\frak m}$. We thus have the following:

\proclaim{Corollary + Definition (3.1.3)} Let $c({\frak m})$
(resp. $c({\frak m}')$) be the dimension 
of the space
$\{ f \in M_k(\varGamma_1(N);{\frak k})_{\frak m} \mid
{\text  {$f$ {\rm has a companion}}}\}$ (resp. $\{ g \in 
M_{k'}(\varGamma_1(N);{\frak k})_{{\frak m}'} \mid
{\text {$g$ {\rm has a companion}}}\}$)
over ${\frak k}$ . Then $c({\frak m})\leq
c({\frak m}')$, and the equality holds when
$k\not= p-1$. 
\endproclaim
Note that we always have $c({\frak m}),~
c({\frak m}') \geq 1$ by (3.1.1).

\proclaim{Theorem (3.1.4)} Let the notation and the assumption be
as above, and suppose that ${\frak r}$ contains a primitive
$N$-th root of unity. 
If $c({\frak m}')=1$, then there is an isomorphism:
$$M_k(\varGamma_1(N);{\frak k})_{\frak m}/
{\frak k}E_k(\chi ,\psi )\overset\sim\to\rightarrow
{\roman{Hom}}_{\frak k}(S_k(\varGamma_1(N);{\frak k})_{\frak m},{\frak k})$$  
of $H_k(\varGamma_1(N);{\frak k})_{\frak m}$-modules.
\endproclaim
\demo{Proof} Our assumption implies that $c({\frak m})=1$, and hence
the twisted version of Ulmer's
pairing ((2.5.1), (2.5.5)) gives us the desired isomorphism.
\qed
\enddemo

\vskip3mm\flushpar
{\bf 3.2. Numerical criterion.} We note that 
an isomorphism stated in (3.1.4) exists only when $p \nmid
\varphi (N)$, for otherwise there is an Eisenstein series different
from, but congruent (mod $\varpi$) to $E_k(\chi ,\psi )$, so that 
the two spaces above cannot have the same dimension.

So, from now on, until the end of this paper, we assume that
{\it{$p$ does not divide $\varphi (N)$}}.
In this subsection, we give a sufficient
condition for $c({\frak m}')$ (actually, $\dim_{\frak k}
M_{k'}(\varGamma_1(N);
{\frak k})_{{\frak m}'}$ itself) to be one, as an application of our
previous work [O4].

We keep the notation of the previous subsection, and set
$$\cases d':=k'-2; \\
\vartheta := \chi \omega^d, \enskip
\vartheta^{\prime}:= \psi \omega^{d'} ; \\
{\frak M}:={\frak M}(\vartheta, \psi ), \enskip
{\frak M}':={\frak M}(\vartheta^{\prime}, \chi )
\enskip {\text{(cf. (1.2.9))}}. 
\endcases $$
Our convention at the beginning of 3.1 excludes 
the case where $k=p-1$ and $\chi = \psi = {\pmb{1}}$
(in which case $k'=2$ and $M_2(SL_2({\bold Z});{\frak r})
=\{ 0\}$).
However, we remark that
the argument below, up to the last paragraph 
in the proof of (3.2.3), works equally in this case.

Now $0 \leq d' \leq p-4$, and the pair $(\vartheta^{\prime},
\chi )$ is not exceptional in the sense of [O4], (1.4.10).
There is a canonical exact sequence:
$$0 \to e\, S(N;\Lambda_{\frak r})_{{\frak M}'} \to 
e\, M(N;\Lambda_{\frak r})_{{\frak M}'} \to
\Lambda_{\frak r} \to 0 \tag 3.2.1$$
which splits uniquely as modules
over $e\,{\Cal H}(N;{\frak r})_{{\frak M}'}$ when tensored with the quotient
field of $\Lambda_{\frak r}$ over $\Lambda_{\frak r}$
(cf. [O4], (1.5.5)). 

\proclaim{Lemma (3.2.2)} The congruence module attached to
{\rm{(3.2.1)}} (cf. {\rm{[O4], 1.1}})
vanishes if and only if 
$e\, S(N;\Lambda_{\frak r})_{{\frak M}'}
=\{ 0\}$, i.e. rank${}_{\Lambda_{\frak r}}\,e\, 
M(N;\Lambda_{\frak r})_{{\frak M}'}=1$.
\endproclaim
\demo{Proof} The $``$if" part is clear. On the other hand, if we assume the
vanishing of the congruence module, the sequence (3.2.1) splits
as modules over $e\, {\Cal H}(N;{\frak r})_{{\frak M}'}$. 
Hence the sequence
tensored with $\Lambda_{\frak r}/(T,\varpi )={\frak k}$ over
$\Lambda_{\frak r}$:
$$0 \to e\, S_2(\varGamma_1(Np);{\frak k})_{{\frak M}'} \to 
e\, M_2(\varGamma_1(Np);{\frak k})_{{\frak M}'} \to {\frak k} \to 0 $$
also splits as modules over the Hecke algebra. The splitting image
of $1 \in {\frak k}$ is a non-zero constant multiple of the unique
Eisenstein series (mod ${\varpi}$) in $e\, M_2(\varGamma_1(Np);\!
{\frak k})_{{\frak M}'}$, which is not a cusp form. 

Suppose that there is a non-zero element $f \in e\, S_2(\varGamma_1(Np);
{\frak k})_{{\frak M}'}$.
Let $\overline{\frak m}'$ be the maximal ideal of $e\, H_2(\varGamma_1(Np);
{\frak k})_{{\frak M}'}$, i.e. the image of ${\frak M}'$ in it, and choose
a non-negative integer $t$ such that $\overline{\frak m}^{\prime t+1}f
=\{ 0\}$ but $\overline{\frak m}^{\prime t}f\not= \{ 0\}$. Then
any non-zero element of $\overline{\frak m}^{\prime t}f$ 
is a cusp form sharing the same eigenvalues for all Hecke
operators as the above Eisenstein series. It follows that a non-zero
constant must belong to $e\, M_2(\varGamma_1(Np);{\frak k})_{{\frak M}'}$.
In view of (1.3.5) and its proof, a non-zero constant also belongs to
$M_{k'}(\varGamma_1(N);{\frak k})$ (resp. $M_{p+1}
(\varGamma_1(N);{\frak k})$) when $d'\not=0$ (resp. $d'=0$). This is 
impossible since $p>3$.
\qed \enddemo

\proclaim{Theorem (3.2.3)} Let the notation be as above. If
the $p$-adic integer
$$\left( \underset{l \nmid{\roman{cond}}( \chi^{-1} \psi)}
\to  {\prod_{l \mid N}} (l^{k'}-(\chi \psi^{-1} )(l))
\right)B_{k', \chi^{-1} \psi} \in {\frak r}$$
is a unit, then $\dim_{\frak k}M_{k'}(\varGamma_1(N);{\frak k})_{{\frak m}'}
=1$, and hence especially 
$c({\frak m}^{\prime})=1$.
Here, $B_{k', \chi^{-1}\psi }$ is the 
generalized Bernoulli number.
\endproclaim
\demo{Proof} By [O4], (1.5.5), the congruence module
attached to the exact sequence (3.2.1) is isomorphic to
$\Lambda_{\frak r}/(A(T; \vartheta^{\prime}, \chi ))$ with
$$A(T;\vartheta^{\prime},\chi )=
\left(\underset{l \nmid {\roman{cond}}(\vartheta^{\prime} \chi^{-1})}\to
{\prod_{l \mid N}} (l^{-1}A_l(T)-(\vartheta^{\prime}\chi^{-1})^{-1}
(l) l^{-2})\right)G(T,\vartheta^{\prime}\chi^{-1}\omega^2) 
\in \Lambda_{\frak r}$$
where $A_l(T)$ and $G(T,\vartheta^{\prime}\chi^{-1}\omega^2 )$ are
given by (1.2.5) and (1.2.6), respectively. 
The value stated in the theorem is a unit multiple of
$A(\gamma^{d'}-1; \vartheta^{\prime} ,\chi )$, and 
hence our assumption implies that
rank$_{\Lambda_{\frak r}}e\, M(N; \Lambda_{\frak r})_{{\frak M}'}=1$
by the previous lemma.

Thus $e\,M(N;\Lambda_{\frak r})_{{\frak M}'}/\omega_{d'}
\cong e\, M_{k'}(\varGamma_1(Np); {\frak r})_{{\frak M}'}$
is a free ${\frak r}$-module of rank one. 
It follows that $M_{k'}(
\varGamma_1(N);{\frak k})_{{\frak m}^{\prime}}$, being a subspace of
$e\,M_{k'}(\varGamma_1(Np);{\frak k})_{{\frak M}^{\prime}}$, is
also one-dimensional.
\qed \enddemo

Note that, when $d'=0$ and $\chi = \psi = {\pmb 1}$, 
the argument above together with
the fact that $B_2=1/6$ assures us that 
$e\, {\Cal H}(1;{\frak r})_{{\frak M}({\pmb 1},{\pmb 1})}=
\Lambda_{\frak r}$.

\vskip3mm
\flushpar
{\bf 3.3. Applications to the structure of $p$-adic Hecke algebras.}
We keep the notation and the assumptions in the previous subsections,
and discuss the Gorenstein property of the Eisenstein components
of the Hecke algebras $e\, {\Cal H}(N;{\frak r})_{\frak M}$ and
$e\, h(N;{\frak r})_{\frak M}$. 
For basic facts about Gorenstein rings, see Bruns and Herzog [BH] or
Eisenbud [E].

Now a local Noetherian ring is a Gorenstein ring if and
only if its quotient by an ideal generated by a regular sequence
is Gorenstein.
So, it follows from (1.5.2) and (1.5.4) that
$e\, {\Cal H}(N;{\frak r})_{\frak M}$ is a Gorenstein ring
if and only if $H_k(\varGamma_1(N);{\frak k})_{\frak m}$ is.
It also follows from (1.5.3) and (1.5.4) that
$M_k(\varGamma_1(N);
{\frak k})_{\frak m}$ is isomorphic to Hom${}_{\frak k}
(H_k(\varGamma_1(N);{\frak k})_{\frak m},{\frak k})$
as a module over $H_k(\varGamma_1(N);{\frak k})_{\frak m}$,
i.e. it is isomorphic to the canonical module of
$H_k(\varGamma_1(N);{\frak k})_{\frak m}.$ Thus 
$e\, {\Cal H}(N;{\frak r})_{\frak M}$ is Gorenstein
if and only if
$M_k(\varGamma_1(N);{\frak k})_{\frak m}$ is a free module
of rank one over $H_k(\varGamma_1(N);{\frak k})_{\frak m}.$

Similarly, when $e\, h(N;{\frak r})_{\frak M}$ is not a
zero-ring, it is Gorenstein if and only
if $S_k(\varGamma_1(N);{\frak k})_{\frak m}$ is free of rank one
over $h_k(\varGamma_1(N);{\frak k})_{\frak m}$.

As before, we use the same symbol $E_k(\chi ,\psi )$
to denote its image in $M_k(\varGamma_1(N);{\frak k})$, in the
following 
\proclaim {Lemma (3.3.1)} The notation being as above,
$H_k(\varGamma_1(N);{\frak k})_{\frak m}$ 
is a Gorenstein ring if
and only if $M_k(\varGamma_1(N);{\frak k})_{\frak m}/
{\frak k}E_k(\chi ,\psi )$ is a cyclic
$H_k(\varGamma_1(N);{\frak k})_{\frak m}$-module.
\endproclaim
\demo{Proof} We only need to prove the $``$if" part; and
assume that the latter condition is satisfied.
If the module in question reduces to $\{ 0 \}$, $H_k(\varGamma_1(N);
{\frak k})_{\frak m}={\frak k}$ is clearly Gorenstein. 
Otherwise, there is a non-zero $f \in M_k(\varGamma_1(N);
{\frak k})_{\frak m}$ such that
$$M_k(\varGamma_1(N);{\frak k})_{\frak m} =
{\frak k}E_k(\chi ,\psi )+ H_k(\varGamma_1(N);{\frak k})_{\frak m}f.$$
Take a non-negative integer $t$ so that ${\frak m}^{t+1}f=\{ 0 \}$
but ${\frak m}^tf \not= \{ 0 \}$. Then any non-zero element of
${\frak m}^tf$ is annihilated by ${\frak m}$. 
Since $M_k(\varGamma_1(N);{\frak k})_{\frak m}$ does not
contain a non-zero constant by (1.5.3), 
such an element must be a
constant multiple of $E_k(\chi , \psi )$. This shows that $f$
generates $M_k(\varGamma_1(N);{\frak k})_{\frak m}$ over
$H_k(\varGamma_1(N);{\frak k})_{\frak m}$.
\qed \enddemo

\proclaim{Theorem (3.3.2)} If $c({\frak m}')=1$, then $e\, {\Cal H}
(N;{\frak r})_{\frak M}$ is a Gorenstein ring. Especially, if the
numerical condition in $(3.2.3)$ is satisfied, $e\, {\Cal H}
(N;{\frak r})_{\frak M}$ is Gorenstein.
\endproclaim
\demo{Proof} By the well-known duality between $S_k(\varGamma_1(N);{\frak k}
)_{\frak m}$ and $h_k(\varGamma_1(N);{\frak k})_{\frak m}$, 

\flushpar
${\roman{Hom}}_{\frak k}(S_k(\varGamma_1(N);{\frak k})_{\frak m},{\frak k})$
is a cyclic module over $H_k(\varGamma_1(N);{\frak k})_{\frak m}$.
So, our conclusion follows from (3.1.4) and (3.3.1).
\qed \enddemo

\demo{Remark {\rm{(3.3.3)}}}
Here is a remark on the relation with a work of Skinner and Wiles
[SW]. We continue to assume that $p \nmid \varphi (N)$ and moreover
that $\psi = {\pmb 1}$. The theorem above implies
 that if $B_{k', \chi^{-1}}\in {\frak r}^{\times}$, then
$e\, {\Cal H}(N;{\frak r})_{\frak M}$, or equivalently
$e\, H_2(\varGamma_1(Np);{\frak r})_{\frak M}$, is Gorenstein. 
We note that Skinner and Wiles have already shown, among others,
that the same
hypothesis implies that $e\, H_2(\varGamma_1(Np);{\frak r})_{\frak M}$,
and hence $e\, {\Cal H}(N;{\frak r})_{\frak M}$ also,
is {\it even a complete intersection}. The reason is as follows:

Let $\Sigma$ be the set of all primes dividing $Np$, and take
$\chi \omega^{d+1}$ as the character $\widetilde{\chi}$ in [SW].
Then the condition $B_{k',\chi^{-1}}\in {\frak r}^{\times}$ is 
equivalent to the vanishing of the
$\widetilde{\chi}$-eigenspace of the $p$-part of the ideal
class group of the splitting field of $\widetilde{\chi}$. This
is one of the basic assumptions in [SW] to set up the deformation problem.
Namely, let ${\bold Q}_{\Sigma}$ be the maximal extension of 
${\bold Q}$ unramified outside $\Sigma$ and the infinite prime.
Skinner and Wiles started with the indecomposable residual 
representation $\rho_0$
of Gal$({\bold Q}_{\Sigma}/{\bold Q})$ of the form
$$\rho_0 =\bmatrix \overline{\widetilde{\chi}} & * \\ 0 & 1 \endbmatrix$$
($\overline{\widetilde{\chi}}$ being the reduction of $\widetilde{\chi}$
modulo $\varpi$) over ${\frak k}$, 
which is unique up to equivalence under the assumption above.
They then considered various types of 
${\frak r}$-deformations of $\rho_0$. 

The Hecke ring ${\bold T}_{\Sigma ,{\frak r}}^{min}$ was defined
as the universal algebra controlling {\it modular} $\Sigma$-minimal
deformations of $\rho_0$.
Concretely, under the assumptions above,
it is the ${\frak r}$-subalgebra of 
$H_2(\varGamma_1(Np);{\frak r})$ generated by
$T(l)$ with prime numbers $l \nmid Np$ and the diamond operators,
localized at its maximal ideal 
corresponding to $E_2(\chi \omega^d,
{\pmb 1})$ (cf. [SW], page 10523).
It was shown in [SW], page 10526, that this ring
is a complete intersection, in their course of proving that
${\bold T}_{\Sigma ,{\frak r}}^{min}$ in fact coincides with
the universal $\Sigma$-minimal deformation ring for $\rho_0$.

Let ${\Cal F}$ be the set of primitive cusp forms belonging to
${\bold T}_{\Sigma, {\frak r}}^{min}$, i.e. those
corresponding to the minimal prime ideals of
${\bold T}_{\Sigma, {\frak r}}^{min}$, and
$\rho_f :{\roman{Gal}}({\bold Q}_{\Sigma}/{\bold Q})\to
GL_2( A_f)$ the $p$-adic representation attached to $f \in {\Cal F}$.
For each
$f \in {\Cal F}$, the $n$-th Fourier coefficient $a(n;f)$ of $f$
is congruent to that of $E_2(\chi \omega^d, {\pmb 1})$ modulo
the maximal ideal of $A_f$ at least
when $n$ is prime to $Np$.
But by applying results of Langlands and Wiles 
to $\rho_f$ as in Wiles
[Wi], pages 500 and 507, one sees that 
$a(q; f)$ for primes $q \mid N$ and $a(p;f)$ are congruent to one modulo
the maximal ideal of $A_f$; i.e. $f$ itself is congruent to
$E_2(\chi \omega^d,{\pmb 1})$. It then
follows that ${\bold T}_{\Sigma ,{\frak r}}^{min}$ coincides with the
${\frak r}$-subalgebra
of $e\, H_2(\varGamma_1(Np);{\frak r})_{\frak M}$ generated by
$T(l)$ with $l \nmid Np$. 

Let $q$ be a prime factor of $N$. Then, again using a result of
Langlands as in [Wi], Remark 2.11,
$T(q)$ belongs to ${\bold T}_{\Sigma ,{\frak r}}^{min}$.
The redundance of $T(p)\in
e\, H_2(\varGamma_1(Np);{\frak r})_{\frak M}$ follows from
(2.5.4) and (1.3.5); or from a result of Wiles as in [Wi], page 507.
Consequently, ${\bold T}_{\Sigma ,{\frak r}}^{min}$ indeed coincides with the
full Hecke ring $e\, H_2(\varGamma_1(Np);{\frak r})_{\frak M}$.

Finally, it follows from this that 
$e\, {\Cal H}(N;{\frak r})_{\frak M}$ is generated over $\Lambda_{\frak r}$
by all $T(l)$ with prime numbers $l$ not dividing $Np$, by
Nakayama's lemma.
\enddemo

\vskip2mm
\demo{Example {\rm{(3.3.4)}}} Assume $N=1$. In this case, there are $(p-1)/2$
$\Lambda$-adic Eisenstein series ${\Cal E}(\omega^d,{\pmb{1}})$
for $0 \leq d \leq p-3$ even. Set ${\frak M}_d:={\frak M}(\omega^d,{\pmb{1}})$.
We have 
$e\, {\Cal H}(1;{\frak r})_{{\frak M}_0}=
\Lambda_{\frak r}$ by the remark at the end of 3.2.
On the other hand, we also have $e\, {\Cal H}(1;{\frak r})_{{\frak M}_{p-3}}= 
\Lambda_{\frak r}$, since $e\, S(1; \Lambda_{\frak r})_{{\frak M}_{p-3}}
= \{ 0 \}$ (cf. [O4], Appendix A.1).

Therefore, so far as the Gorenstein property is concerned, we only need
to treat the case where $2 \leq d \leq p-5$, or $4 \leq k \leq
p-3$. 
For this range of $d$, the roles
of $k$ and $k'$ are symmetric, and the 
numerical criterion above may be 
rephrased as follows: If either one of $B_k$ or $B_{k'}$ is not divisible
by $p$, then both $e\, {\Cal H}(1;{\frak r})_{{\frak M}_d}$ and
$e\, {\Cal H}(1;{\frak r})_{{\frak M}_{d'}}$ are Gorenstein
(actually complete intersections by Skinner and Wiles).
Thus one can naturally ask if there is a pair $(k,k')$
$(4 \leq k \leq p-3, k+k'=p+1)$ such that $B_k\equiv B_{k'}
\equiv 0$ (mod $p$). As is well-known, when $p \equiv 3$ (mod 4),
$B_{(p+1)/2} \not\equiv 0$ (mod $p$) (cf. Washington [Wa], Exercise
5.9); and it can be easily checked that there is no such a pair
for $p \leq 4001$, using the table in [Wa], Tables, Section 2. 
At present, we know no example such that $B_k\equiv B_{k'}\equiv 0$ 
(mod $p$), and hence neither know whether or not the implication:
$$B_k \not\equiv 0 \enskip ({\roman{mod}} ~~p) \enskip {\roman{or}}\enskip
B_{k'}\not\equiv 0 \enskip
({\roman{mod}} ~~p)\enskip \Rightarrow \enskip 
c({\frak m})=c({\frak m}')=1 $$
is strict.
\enddemo

\vskip2mm
We next consider the structure of the cuspidal Hecke algebra
$e\, h(N;{\frak r})_{\frak M}$ when $c({\frak m}')=1$.
It is based on the isomorphism:
$$S_k(\varGamma_1(N);{\frak k})_{\frak m}\cong
{\roman{Hom}}_{\frak k}(M_k(\varGamma_1(N);{\frak k})_{\frak m}
/{\frak k}E_k(\chi ,\psi ),{\frak k})
\tag 3.3.5$$
deduced from (3.1.4).

\proclaim{Lemma (3.3.6)} Let $\overline{\frak m}$ be the image
of ${\frak m}$ in $H_k(\varGamma_1(N);{\frak k})_{\frak m}$. 
If $c({\frak m}')=1$, $S_k(\varGamma_1(N);{\frak k})_{\frak m}$
is isomorphic to $\overline{\frak m}$ as an $H_k(\varGamma_1(N);
{\frak k})_{\frak m}$-module.
\endproclaim
\demo{Proof} Taking the ${\frak k}$-dual of the exact sequence:
$$0 \to {\frak k}E_k(\chi ,\psi ) \to M_k(\varGamma_1(N);{\frak k})_{\frak m}
\to M_k(\varGamma_1(N);{\frak k})_{\frak m}/{\frak k}E_k(\chi ,\psi)
\to 0$$
we see that ${\roman{Hom}}_{\frak k}(M_k(\varGamma_1(N);{\frak k})_{\frak m}
/{\frak k}E_k(\chi ,\psi ),{\frak k})$ is isomorphic to the
kernel of the homomorphism $H_k(\varGamma_1(N);{\frak k})_{\frak m}
\to {\frak k}$ sending $t$ to $a(1;E_k(\chi ,\psi )\mid t)$,
by (1.5.3) and (1.5.4).
Our claim follows from the isomorphism above.
\qed \enddemo

\proclaim{Corollary (3.3.7)} Let $\overline{\frak m}_0$ be the
image of ${\Cal I}(\vartheta ,\psi )_{\frak M}$ (cf. $(1.2.9)$) 
in 
\flushpar
$e\, H_k(\varGamma_1(N)\cap
\varGamma_0(p);{\frak k})_{\frak M}$ via $(1.5.2)$.
If $c({\frak m}')=1$, then $e\, S_k(\varGamma_1(N)\cap
\varGamma_0(p);{\frak k})_{\frak M}$ is isomorphic to $\overline{\frak m}_0$
as an $e\, H_k(\varGamma_1(N)\cap \varGamma_0(p);{\frak k})_{\frak M}$-module.
\endproclaim
\demo{Proof} We know that $M_k(\varGamma_1(N);{\frak k})_{\frak m}
=e\, M_k(\varGamma_1(N)\cap \varGamma_0(p);{\frak k})_{\frak M}$,
as noted in the proof of (1.5.3), and the same holds for cusp forms. Also,
we have a canonical isomorphism: 
$H_k(\varGamma_1(N);{\frak k})_{\frak m}\cong
e\, H_k(\varGamma_1(N)\cap \varGamma_0(p);{\frak k})_{\frak M}$
by (1.5.4). Our claim follows immediately from (3.3.6).
\qed \enddemo

\proclaim{Theorem (3.3.8)} Assume that $e\,h(N;{\frak r})_{\frak M}$
is not a zero-ring.
If $c({\frak m}')=1$, then
the following conditions are equivalent:

{\rm{i)}}~$e\, h(N;{\frak r})_{\frak M}$ is a Gorenstein ring;

{\rm{ii)}}~$e\, h(N;{\frak r})_{\frak M}$ is a complete intersection;

{\rm{iii)}}~the Eisenstein ideal ${\Cal I}(\vartheta ,\psi )_{\frak M}$ of
$e\, {\Cal H}(N;{\frak r})_{\frak M}$ is principal;

{\rm{iv)}}~the image $I$ of ${\Cal I}(\vartheta ,\psi )_{\frak M}$ 
in $e\, h(N;{\frak r})_{\frak M}$, the Eisenstein ideal of $e\, h(N;
{\frak r})_{\frak M}$,
is principal.
\endproclaim
\demo{Proof} The kernel of the natural surjection: $e\, {\Cal H}(N;
{\frak r})_{\frak M}\to e\, h(N;{\frak r})_{\frak M}$ has trivial
intersection with ${\Cal I}(\vartheta ,\psi )_{\frak M}$ because an element
in the intersection annihilates both ($\Lambda$-adic) cusp forms and the
Eisenstein series. Thus the conditions iii) and iv) are equivalent.
That iv) $\Rightarrow$ ii) is clear, and that ii) $\Rightarrow$ i) is 
well-known.
Note that these implications hold unconditionally.

We now show that the condition i) implies iii). 
The condition i) is equivalent to the cyclicity of the 
$e\, H_k(\varGamma_1(N)\cap \varGamma_0(p);{\frak k})_{\frak M}$-module
$e\, S_k(\varGamma_1(N)\cap \varGamma_0(p);{\frak k})_{\frak M}$.
By (3.3.7), this is in turn equivalent to the principality of the ideal
$\overline{\frak m}_0$ of

\flushpar
$e\, H_k(\varGamma_1(N)\cap \varGamma_0(p);
{\frak k})_{\frak M}$. Consider the following diagram of Hecke
algebras:
$$\CD
e\, {\Cal H}(N;{\frak r})_{\frak M}@>>>
e\, {\Cal H}(N;{\frak r})_{\frak M}/\omega_d @>>>
e\, {\Cal H}(N;{\frak r})_{\frak M}/(\omega_d, \varpi )\\
@. @V(1.5.2)V{\wr}V @VV{\wr}V \\
@. e\,H_k(\varGamma_1(N)\cap \varGamma_0(p);{\frak r})_{\frak M}
@>>> e\, H_k(\varGamma_1(N)\cap \varGamma_0(p);{\frak k})_{\frak M}.
\endCD $$
Let $I_k$ be the image of ${\Cal I}(\vartheta ,\psi )_{\frak M}$ 
in $e\, H_k(\varGamma_1(N)
\cap \varGamma_0(p);{\frak r})_{\frak M}$. It is the annihilator of
$E_k(\chi_1 ,\psi ) \in e\, M_k(\varGamma_1(N)\cap \varGamma_0(p);
{\frak r})_{\frak M}$, and hence $I_k\cap \varpi e\, H_k(\varGamma_1(N)
\cap \varGamma_0(p);{\frak r})_{\frak M}=\varpi I_k$. It follows that
$I_k/\varpi I_k$ is isomorphic to $\overline{\frak m}_0 $ as an
$e\, H_k(\varGamma_1(N)\cap \varGamma_0(p);{\frak r})_{\frak M}$-module.
Thus Nakayama's lemma implies that $I_k$ is principal if the
condition i) holds. The same reasoning applies to the left situation
in the diagram above, showing that ${\Cal I}(\vartheta ,\psi )_{\frak M}$ is
also principal under i).
\qed
\enddemo

We note that $e\,h(N;{\frak r})_{\frak M}\not= \{ 0\}$ if and only if
the quotient of $e\, h(N;{\frak r})_{\frak M}$ by
the image of ${\Cal I}(\vartheta ,\psi )_{\frak M}$
is non-zero. By [O4], 3.2 and (1.5.5),
this is the case exactly when $G(T,\vartheta \omega^2)\not\in
\Lambda_{\frak r}^{\times}$.

\vskip3mm \flushpar
{\bf 3.4. Relation with the theory of cyclotomic
fields.} In this final subsection, we throughout assume that
$\psi ={\pmb 1}$. Thus $\chi$ is a primitive Dirichlet character
of conductor $N$ such that $\chi (-1)= (-1)^d$,
and $\vartheta = \chi \omega^d$ with
$1 \leq d \leq p-3$.
We take ${\frak r}$ to be the ring generated
by the values of $\chi$ over ${\bold Z}_p$.
We are going to consider consequences of the Gorenstein
property of $e\, {\Cal H}(N;{\frak r})_{\frak M}$
on the two-dimensional
Galois representation attached to $e\, S(N; \Lambda_{\frak r})_{\frak M}$.

To do this, we use the same terminology as in [O2] and 
[O3]. Thus $e^*ES_p(N)_{\frak r}$ (resp. $e^*GES_p
(N)_{\frak r}$) stands for the ordinary part of
$\underset{r\geq 1}\to\varprojlim \,
H_{\roman{\acute{e}t}}^1(X_1(Np^r)\otimes_{\bold Q}\overline{\bold Q},
{\bold Z}_p)\otimes_{{\bold Z}_p}{\frak r}$ 
(resp. $\underset{r\geq 1}\to\varprojlim \,
H_{\roman{\acute{e}t}}^1(Y_1(Np^r)\otimes_{\bold Q}\overline{\bold Q},
{\bold Z}_p)\otimes_{{\bold Z}_p}{\frak r}$), 
and $e^*{\Cal H}^*(N;{\frak r})$ (resp. $e^*
h^*(N;{\frak r})$) for Hida's universal ordinary $p$-adic Hecke algebra
acting on this space. This algebra is canonically isomorphic to
$e\,{\Cal H}(N;{\frak r})$ (resp. $e\, h(N;{\frak r})$). Let
${\frak M}^*$ be the maximal ideal of $e^*{\Cal H}^*(N;
{\frak r})$ corresponding to ${\frak M}$ via this isomophism,
and set:
$$\cases
{\frak h}^* :=e^*h^*(N;{\frak r})_{{\frak M}^*} \\
{\frak H}^*:=e^*{\Cal H}^*(N;{\frak r})_{{\frak M}^*} \\
X:= e^*ES_p(N)_{{\frak r}, {\frak M}^*} \\
Y:= e^*GES_p(N)_{{\frak r},{\frak M}^*}
\endcases \tag 3.4.1$$
for simplicity. In what follows, we will always assume that
$G(T, \vartheta \omega^2)\not\in \Lambda_{\frak r}^{\times}$,
so that ${\frak h}^*$ is not a zero-ring.
Let ${\Cal I}^*$ be the Eisenstein ideal of ${\frak H}^*$,
i.e. the ideal corresponding to ${\Cal I}(\vartheta ,{\pmb 1})_{\frak M}$
via ${\frak H}^*\cong e\,{\Cal H}(N;{\frak r})_{\frak M}$,
and $I^*$ its image in ${\frak h}^*$. We therefore have an isomorphism:
$${\frak h}^*/I^*\cong \Lambda_{\frak r}/(G(T;\vartheta \omega^2)).
\tag 3.4.2 $$

Let $I_p \subset {\roman{Gal}}(\overline{\bold Q}_p/{\bold Q}_p)$
be the inertia group. We know that
$$X^{I_p}= Y^{I_p}=: X_+ \tag 3.4.3$$
is a free ${\frak h}^*$-module of rank one,
and that $ X/X_+$ is isomorphic to the
$\Lambda_{\frak r}$-dual of ${\frak h}^*$ as an ${\frak h}^*$-module
([O3], (2.3.6)).

\proclaim{Lemma (3.4.4)}
$Y/X_+$ is isomorphic to
the $\Lambda_{\frak r}$-dual of ${\frak H}^*$ as an 
${\frak H}^*$-module. \endproclaim
\demo{Proof} We use the same notation as in [O3]. Thus ${\frak o}$
stands for the ring of integers of a sufficiently large complete
subfield of ${\bold C}_p$, and we remind of us that $\Lambda_{\frak o}$
is faithfully flat over $\Lambda_{\frak r}$ ([O3], (2.1.1)).
By [O3], (2.1.11) and (2.1.12), we know that $(Y/X_+)\otimes_
{\Lambda_{\frak r}}\Lambda_{\frak o}$ is isomorphic to 
$e\, M(N; \Lambda_{\frak o})_{\frak M}$. 
It follows from the duality (1.5.3) that the $\Lambda_{\frak o}$-dual 
of $(Y/X_+)\otimes_{{\Lambda}_{\frak r}}\Lambda_{\frak o}$ is free
of rank one over ${\frak H}^*\otimes_{\Lambda_{\frak r}}\Lambda_{\frak o}$.
From the remark above, we conclude that the $\Lambda_{\frak r}$-dual
of $Y/X_+$ itself is free of rank one over ${\frak H}^*$. \qed
\enddemo

Now assume that ${\frak H}^*$ is Gorenstein
(which is implied by the condition 
$c({\frak m}')=1$ by (3.3.2)). It follows from the lemma above
that $ Y/X_+$ is a free 
${\frak H}^*$-module of rank one. Set
$$
\widetilde{X}:=Y\otimes_{{\frak H}^*}{\frak h}^*.
\tag 3.4.5$$
Then we have a commutative diagram of ${\frak h}^*$-modules:
$$\CD
0 @>>> X_+ @>>> X @>>> X/X_+@>>> 0 \quad ({\roman{exact}})\\
@. @| @VVV @VVV @. \\
0 @>>> X_+ @>>> \widetilde{X} @>>>
(Y/X_+)\otimes_{{\frak H}^*}{\frak h}^*
 @>>> 0 \quad ({\roman{split}}).
\endCD \tag 3.4.6$$
Since $Y/X$ is isomorphic to ${\frak H}^*/{\Cal I}^*$ as an
${\frak H}^*$-module by [O2], (5.2.11), the cokernel of the
middle vertical arrow is isomorphic to ${\frak h}^*/I^*$ as an
${\frak h}^*$-module, which is $\Lambda_{\frak r}$-torsion by (3.4.2).
Since $X$ and $\widetilde{X}$ are free $\Lambda_{\frak r}$-modules of the
same rank, the two vertical arrows above are injective.

Fix a $\sigma_0 \in I_p$ which acts as multiplication by 
$\omega^{-d-1}(\sigma_0)\not= 1$ on $\widetilde{X}/X_+$, 
and use it to split the two horizontal sequences in (3.4.6)
as ${\frak h}^*$-modules, as in [O2], 5.3. Especially, if we denote by
$\widetilde{X}_-$ the $\omega^{-d-1}(\sigma_0)$-eigenspace 
with respect to the action of
$\sigma_0$ on $\widetilde{X}$, 
$\widetilde{X}$ is a direct sum:
$$\widetilde{X}=
\widetilde{X}_- \oplus X_+ \tag 3.4.7 $$
of two free ${\frak h}^*$-modules of rank one. Fixing bases of
$\widetilde{X}_-$ and $X_+$ (in this order) respectively, we obtain a
representation
$$\widetilde{\rho}: {\roman{Gal}}(\overline{\bold Q}/{\bold Q})
\to GL_{{\frak h}^*}(\widetilde{X})\cong GL_2({\frak h}^*)\quad
{\roman{with}}\quad \widetilde{\rho}(\sigma )=
\bmatrix \widetilde{a}(\sigma ) &\widetilde{b}(\sigma )\\
\widetilde{c}(\sigma ) & \widetilde{d}(\sigma )
\endbmatrix . \tag 3.4.8$$
From the remark above, this realizes the representation into
$GL_2({\frak h}^*\otimes_{\Lambda_{\frak r}}
{\roman Q}(\Lambda_{\frak r}))$, considered in 
[O2], 5.3, over ${\frak h}^*$,
where ${\roman Q}(\Lambda_{\frak r})$ denotes the
quotient field of $\Lambda_{\frak r}$.

We can then apply the method of Harder and Pink [HP] and
Kurihara [Ku] to this representation as exposed 
in loc. cit.: Let $\widetilde{B}$ (resp. $\widetilde{C}$) be the
ideal of ${\frak h}^*$ generated by $\widetilde{b}(\sigma )$
(resp. $\widetilde{c}(\sigma )$) for all $\sigma \in
{\roman{Gal}}(\overline{\bold Q}/{\bold Q})$. 
By the final remark in (3.3.3), $I^*$
is generated by all $T^*(l) -1- lT^*(l,l)$
with prime numbers $l$ not dividing $Np$, and hence the ideal
denoted by $J^*$ in [O2] coincides with $I^*$. Thus we have
$$\widetilde{B}\widetilde{C}= I^*. \tag 3.4.9 $$
We moreover have:
$$\widetilde{B}=I^* \enskip {\roman{and}} \enskip
\widetilde{C}={\frak h}^*.\tag 3.4.10 $$
Indeed, it is enough to show that $\widetilde{B}$ is
contained in $I^*$. But from the construction of $\widetilde{\rho}$,
every $\widetilde{b}(\sigma )$ annihilates $\widetilde{X}/X
\cong {\frak h}^*/I^*$.

Let $F$ be the abelian extension of ${\bold Q}$ corresponding to
$\vartheta \omega$, and $F_{\infty}$ its cyclotomic
${\bold Z}_p$-extension. Let $L_{\infty}$ be the maximal abelian
pro-$p$ unramified extension of $F_{\infty}$. 
The group $\Delta :={\roman{Gal}}(F/{\bold Q})$ acts on 
${\roman{Gal}}(L_{\infty}/F_{\infty})$ in the usual manner,
and we set
$${\roman{Gal}}(L_{\infty}/F_{\infty})_{(\vartheta \omega )^{-1}}
:={\roman{Gal}}(L_{\infty}/F_{\infty})\otimes_{{\bold Z}_p[\Delta ]}
{\frak r} \tag 3.4.11$$
via the homomorphism ${\bold Z}_p[\Delta ]\to {\frak r}$ induced
by $(\vartheta \omega )^{-1}$. This is a module
over ${\frak r}[[{\roman {Gal}}(F_{\infty}/ F)]]$ in a natural manner.

Recall that from the outset we have fixed a topological generator
$\gamma$ of the multiplicative group $1+p{\bold Z}_p$, and
used it to identify $\Lambda_{\frak r} ={\frak r}[[1+p{\bold Z}_p]]$
with ${\frak r}[[T]]$ via $\gamma \leftrightarrow 1+T$, and hence
a $\Lambda_{\frak r}$-module with a ${\frak r}[[T]]$-module.
Now for a $\Lambda_{\frak r}$-module $M$, we denote by $M^{\dagger}$
the $\Lambda_{\frak r}$-module for which the action of $T$ is
newly given by that of
$\gamma^{-1}(1+T)^{-1}-1$. One easily checks that
this $\Lambda_{\frak r}$-module structure 
of $M^{\dagger}$ does not depend on
the choice of $\gamma$.
On the other hand, the $p$-cyclotomic character
induces an isomorphism:
${\frak r}[[{\roman{Gal}}(F_{\infty}/F)]]\overset{\sim}\to \rightarrow
{\frak r}
[[1+p{\bold Z}_p]]=\Lambda_{\frak r}$. We consider
${\roman{Gal}}(L_{\infty}/F_{\infty})_{(\vartheta \omega )^{-1}}$
as a $\Lambda_{\frak r}$-module through it. With these terminologies,
the following theorem follows
from [O2], (5.3.18) and (5.3.20):

\proclaim{Theorem (3.4.12)} Assume that 
$e\,{\Cal H}(N;{\frak r})_{\frak M}$ is a Gorenstein ring. Then
the Iwasawa module ${\roman{Gal}}
(L_{\infty}/F_{\infty})_{(\vartheta \omega )^{-1}}$ is isomorphic
to $(I^*/I^{*2})^{\dagger}$ as a $\Lambda_{\frak r}$-module. 
\endproclaim

\proclaim{Corollary (3.4.13)} Under the same hypothesis, 
the conditions {\rm{i) - iv)}} in {\rm{(3.3.8)}} are all equivalent to the
cyclicity of the Iwasawa module 
${\roman{Gal}}(L_{\infty}/F_{\infty})_{(\vartheta \omega )^{-1}}$.
\endproclaim 

\demo{Proof} It follows from the theorem above that the cyclicity
of the Iwasawa module is equivalent to the principality of 
the Eisenstein ideal $I$ of $e\, h(N;{\frak r})_{\frak M}$. 
Therefore it is enough to show that the condition i) in (3.3.8) implies
the cyclity.  When $N=1$, this is due to
Harder, Pink and Kurihara,
and their method easily generalizes to
the present case as follows: 
Under that condition, the module $X$ is 
a direct sum of free $ {\frak h}^*$-modules $X_+$ and 
$X_-$ of rank one, $X_-$ being defined in the same manner
as $\widetilde{X}_-$ from $X$. We can then consider the
Galois representation 
$$\rho :{\roman{Gal}}(\overline{\bold Q}/{\bold Q})\to
GL_{{\frak h}^*}(X)\cong GL_2({\frak h}^*) \quad
{\roman{with}} \quad \rho (\sigma ) =  \bmatrix a(\sigma ) &
b(\sigma ) \\ c(\sigma ) & d(\sigma ) \endbmatrix $$
in the same way as (3.4.8).
Let $B$ (resp. $C$) be the ideal of ${\frak h}^*$ generated by
all $b(\sigma )$ (resp. $c(\sigma )$).
For the same reason as (3.4.9), we have $BC=I^*$.

By [O2], (5.2.16), there is a 
${\roman{Gal}}(\overline{\bold Q}/{\bold Q})$-stable 
${\frak h}^*$-submodule $V$ of $X/G(T, \vartheta \omega^2 )$
which is isomorphic to ${\frak h}^*/I^*$. Moreover the element
$\sigma_0 \in I_p$ above acts as multiplication by $\omega^{-d-1}(
\sigma_0)$ on $V$. Therefore $V$ is in fact contained in $X_{-}
/G(T, \vartheta \omega^2)$. 
Since each $c(\sigma )$ annihilates $V$, 
we conclude that $C\subseteq I^*$.

It follws that $B={\frak h}^*$ and $C=I^*$.
Consequently we have:
$${\roman{Gal}}(L_{\infty}/F_{\infty})_{(\vartheta \omega )^{-1}}
\cong (\Lambda_{\frak r}/(G(T, \vartheta \omega^2)))^{\dagger}$$
by [O2], (5.3.18) and (5.3.20). \qed
\enddemo

\vskip5mm
Masami Ohta

Department of Mathematics,

Faculty of Science, Tokai University,

Hiratsuka, Kanagawa, 259-1292, Japan

E-mail: ohta\@sm.u-tokai.ac.jp

\enddocument